\documentclass[twoside,psfig]{article}

\oddsidemargin=\evensidemargin
\catcode`\@=11
\long\def\@makefntext#1{
\protect\noindent \hbox to 3.2pt {\hskip-.9pt  
$^{{\eightrm\@thefnmark}}$\hfil}#1\hfill}               

\def\ps@myheadings{\let\@mkboth\@gobbletwo              
\def\@oddhead{\hbox{}
\rightmark\hfil\eightrm\thepage}   
\def\@oddfoot{}\def\@evenhead{\eightrm\thepage\hfil
\leftmark\hbox{}}\def\@evenfoot{}
\def\sectionmark##1{}\def\subsectionmark##1{}}

\catcode`@=11                                
\def\ps@plain{\let\@mkboth\@gobbletwo
     \def\@oddhead{}\def\@oddfoot{\eightrm\hfil\thepage
     \hfil}\def\@evenhead{}\let\@evenfoot\@oddfoot}

\renewcommand{\thefootnote}{\fnsymbol{footnote}}

\newcounter{sectionc}\newcounter{subsectionc}\newcounter{subsubsectionc}
\renewcommand{\section}[1] {\vspace{12pt}\addtocounter{sectionc}{1} 
\setcounter{subsectionc}{0}\setcounter{subsubsectionc}{0}\noindent 
        {\tenbf\thesectionc. #1}\par\vspace{5pt}}
\renewcommand{\subsection}[1] {\vspace{12pt}\addtocounter{subsectionc}{1} 
        \setcounter{subsubsectionc}{0}\noindent 
        {\bf\thesectionc.\thesubsectionc. 
        {\kern1pt \bfit #1}}\par\vspace{5pt}}
\renewcommand{\subsubsection}[1] {\vspace{12pt}
        \addtocounter{subsubsectionc}{1}
        \noindent
        {\tenrm\thesectionc.\thesubsectionc.\thesubsubsectionc. {\kern1pt 
        \it #1}}\par\vspace{5pt}}
\newcommand{\nonumsection}[1] {\vspace{12pt}\noindent{\tenbf #1}
        \par\vspace{5pt}}

\newcounter{appendixc}
\newcounter{subappendixc}[appendixc]
\newcounter{subsubappendixc}[subappendixc]

\renewcommand{\appendix}[1] {\vspace{12pt}      
        \refstepcounter{appendixc}              
        \setcounter{figure}{0}
        \setcounter{table}{0}
        \setcounter{lemma}{0}
        \setcounter{theorem}{0}
        \setcounter{corollary}{0}
        \setcounter{definition}{0}
        \setcounter{equation}{0}
        \renewcommand{\thefigure}{\Alph{appendixc}.\arabic{figure}}
        \renewcommand{\thetable}{\Alph{appendixc}.\arabic{table}}
        \renewcommand{\theappendixc}{\Alph{appendixc}}
        \renewcommand{\thelemma}{\Alph{appendixc}.\arabic{lemma}}
        \renewcommand{\thetheorem}{\Alph{appendixc}.\arabic{theorem}}
        \renewcommand{\thedefinition}{\Alph{appendixc}.\arabic{definition}}
        \renewcommand{\thecorollary}{\Alph{appendixc}.\arabic{corollary}}
        \renewcommand{\theequation}{\Alph{appendixc}.\arabic{equation}}
        \noindent{\tenbf Appendix \theappendixc #1}\par\vspace{5pt}}

\newcommand{\textlineskip}{\baselineskip=13pt}
\newcommand{\smalllineskip}{\baselineskip=10pt}

\newcommand{\copyrightheading}[1]
        {\vspace*{-2.5cm}\smalllineskip{\flushleft
        {\footnotesize \quad #1}\\ 
        {\footnotesize \quad }\\ 
         }}

\newcommand{\pub}[1]{{\begin{center}\footnotesize\smalllineskip 
        \\ 
        \end{center}
        }}

\def\abstracts#1#2#3#4{{
        \centering{\begin{minipage}{4.5in}\footnotesize\baselineskip=10pt
        \centerline{ABSTRACT} 
        \parindent=15pt #1\par 
        \parindent=15pt #2\par
        \parindent=15pt #3\par
        \parindent=15pt #4\par
        \end{minipage}}\par}} 

\newcounter{itemlistc}
\newcounter{romanlistc}
\newcounter{alphlistc}
\newcounter{arabiclistc}

\newcommand{\fcaption}[1]{
        \refstepcounter{figure}
        \setbox\@tempboxa = \hbox{\footnotesize Fig.~\thefigure. #1}
        \ifdim \wd\@tempboxa > 5in
           {\begin{center}
        \parbox{5in}{\footnotesize\smalllineskip Fig.~\thefigure. #1}
            \end{center}}
        \else
             {\begin{center}
             {\footnotesize Fig.~\thefigure. #1}
              \end{center}}
        \fi}

\newcommand{\tcaption}[1]{
        \refstepcounter{table}
        \setbox\@tempboxa = \hbox{\footnotesize Table~\thetable. #1}
        \ifdim \wd\@tempboxa > 5in
           {\begin{center}
        \parbox{5in}{\footnotesize\smalllineskip Table~\thetable. #1}
            \end{center}}
        \else
             {\begin{center}
             {\footnotesize Table~\thetable. #1}
              \end{center}}
        \fi}

\def\pmb#1{\setbox0=\hbox{#1}
        \kern-.025em\copy0\kern-\wd0
        \kern.05em\copy0\kern-\wd0
        \kern-.025em\raise.0433em\box0}

\def\fnt#1#2{\footnotetext{\kern-.3em
        {$^{\mbox{\scriptsize #1}}$}{#2}}}

\def\fpage#1{\begingroup
\voffset=.3in
\thispagestyle{empty}\begin{table}[b]\centerline{\footnotesize #1}
        \end{table}\endgroup}

\def\runninghead#1#2{\pagestyle{myheadings}
\markboth{{\protect\footnotesize\it{\quad #1}}\hfill}
{\hfill{\protect\footnotesize\it{#2\quad}}}}

\font\tenrm=cmr10
 
\font\tenbf=cmbx10
\font\bfit=cmbxti10 at 10pt
\font\ninerm=cmr9

\font\eightrm=cmr8

\def\@begintheorem#1#2{\trivlist        
        \item[\hskip\labelsep{\bf #1\ #2.}]} 
\def\@opargbegintheorem#1#2#3{\trivlist
        \item[\hskip\labelsep{\bf #1\ #2\ (#3).}]}

        {\setcounter{itemlistc}{0}              
         \begin{list}{$\bullet$}                
        {\usecounter{itemlistc}                 
         \leftmargin10pt               
         \setlength{\parsep}{0pt}
         \setlength{\itemsep}{0pt}     
        }}{\end{list}}

        {\setcounter{romanlistc}{0}             
         \begin{list}{$($\roman{romanlistc}$)$} 
        {\usecounter{romanlistc}                
         \leftmargin18pt 
         \setlength{\parsep}{0pt}
         \setlength{\itemsep}{0pt}      
         \settowidth{\labelwidth}{#1}                          
        }}{\end{list}}

        {\setcounter{enumii}{0}                 
         \begin{list}{$($\alph{enumii}$)$}      
        {\usecounter{enumii}                    
         \leftmargin18pt                
         \setlength{\parsep}{0pt}
         \setlength{\itemsep}{0pt}      
         \settowidth{\labelwidth}{#1}                          
        }}{\end{list}}

\def\qed{\hbox{${\vcenter{\vbox{                        
   \hrule height 0.4pt\hbox{\vrule width 0.4pt height 6pt
   \kern5pt\vrule width 0.4pt}\hrule height 0.4pt}}}$}}

\renewcommand{\thefootnote}{\fnsymbol{footnote}}  

\def\theequation{\thesectionc.\arabic{equation}}  

 \usepackage{amsfonts}
 \usepackage{amssymb}
 \usepackage{amsmath}
 \usepackage[dvips]{epsfig,graphics}
 \usepackage{latexsym, amsfonts, longtable, tabularx}
 \usepackage{epsfig} 
 \usepackage{psfrag}  
 \usepackage{portland}

\newcommand{\dihedral}{\ensuremath{\mathbb{D}_2}}
\newcommand{\cyclic}{\ensuremath{\mathbb{Z}_2}}

\begin{document}
\setlength{\textheight}{7.7truein}  

\runninghead{A. Champanerkar, I. Kofman \& E. Patterson}
{The next simplest hyperbolic knots}

\normalsize\textlineskip
\thispagestyle{empty}
\setcounter{page}{1}

\copyrightheading{}

\vspace*{0.88truein}

\fpage{1}
\centerline{\bf THE NEXT SIMPLEST HYPERBOLIC KNOTS}
\baselineskip=13pt
\vspace*{0.37truein}
\centerline{\footnotesize ABHIJIT CHAMPANERKAR}
\baselineskip=12pt
\centerline{\footnotesize\it Department of Mathematics, Barnard College, Columbia University}
\baselineskip=10pt
\centerline{\footnotesize\it New York, NY 10027}

\vspace*{10pt}
\centerline{\footnotesize ILYA KOFMAN}
\baselineskip=12pt
\centerline{\footnotesize\it Department of Mathematics, Columbia University}
\baselineskip=10pt
\centerline{\footnotesize\it New York, NY 10027}

\vspace*{10pt}
\centerline{\footnotesize ERIC PATTERSON}
\baselineskip=12pt
\centerline{\footnotesize\it Department of Mathematics, University of Chicago}
\baselineskip=10pt
\centerline{\footnotesize\it 5734 S. University Avenue, Chicago, Illinois 60637}

\vspace*{0.225truein}
\ \\ 

\vspace*{0.21truein} 
\abstracts{
We complete the project begun by Callahan, Dean and Weeks to identify all knots whose
complements are in the SnapPea census of hyperbolic manifolds with seven
or fewer tetrahedra.  Many of these ``simple'' hyperbolic knots have high
crossing number.  We also compute their Jones polynomials.
}{}{}{}



\vspace*{1pt}\textlineskip      

\section{Introduction}
\vspace*{-0.5pt}

Computer tabulation of knots, links and their invariants, such as the
Jones polynomial, has led to many advances.
For both practical and historical reasons, knots continue to be identified
according to their crossing number, despite the difficulty to compute it.
On the other hand, computer programs such as SnapPea \cite{snappea}, which have been
essential to investigate volume and other hyperbolic invariants of
3-manifolds, have become important tools to study knots via the geometry
of their complements.
One indication that these seemingly disparate invariants are strongly
related is the remarkable conjecture that the colored Jones polynomials
determine the hyperbolic volume of the knot complement \cite{mm2001}.

It is therefore natural to ask which knots have complements with the
simplest hyperbolic geometry.
In \cite{cdw99}, all hyperbolic knot complements with 6 or fewer tetrahedra
were found, and the corresponding 72 knots were identified.
Of the 4587 orientable  one-cusped hyperbolic manifolds with at most seven ideal
tetrahedra, 3388 have exactly seven tetrahedra, and completing this
project has required additional techniques.

In an undergraduate summer project at Columbia University, Gerald Brant,
Jonathan Levine, Eric Patterson and  Rustam Salari followed
\cite{cdw99} to determine which census manifolds are knot complements.
Using Gromov and Thurston's $2\pi$ theorem, which has been recently
sharpened in \cite{lackenby2000}, there are at most 12 Dehn fillings, 
all on short
filling curves, which can result in $S^3$ from any 1-cusped hyperbolic manifold.

Our computer program, which used SnapPea to simplify the
fundamental group, immediately identified 128 knot complements in the
census.
Another 150 census manifolds had fillings with trivial first homology, but
SnapPea could not simplify their fundamental groups.
We used Testisom \cite{testisom} to show that 148 of these were nontrivial.  Of the
remaining two census manifolds, $v220$ turned out to be a knot complement.
For the final manifold in doubt, $(0,1)$-filling on $v312$, Snap \cite{snap} was used
to compute the base orbifold for a flat structure: it is Seifert fibered
with 3 singular fibers over a sphere, and hence not $S^3$.  This left us
with 129 complements of knots in a homotopy 3-sphere.

To find the corresponding knots in $S^3$, we first searched through the
computer tabulation of knots up to 16 crossings, provided to us in an
electronic file by Morwen Thistlethwaite.
Using this vast tabulation, which is also available via the program
Knotscape \cite{Knotscape}, we found 32 knots.

We then ran programs for several weeks to compute geometric invariants of
{\em twisted torus knots}.
The twisted torus knot $T(p,q,r,s)$ is obtained by performing $s$ full
twists on $r$ parallel strands of a $(p,q)$ torus knot, where $p$ and $q$
are coprime and $1< r < p$ (see \cite{cdw99}). We denote the mirror image of $T(p,q,r,s)$ by $T(p,-q,r,-s)$.
These are easily expressed as closed braids, and we computed their
complements using Nathan Dunfield's python addition to SnapPea for closed
braids \cite{dunfield}.
In this way, we found 72 more knots.

Knots for the remaining 25 complements were obtained as follows.  Snap and
SnapPea have a library of links with up to ten crossings.
We ran a program to compute Dehn fillings on $(n-1)$-component sublinks of
these $n$-component links to find matches with the given knot complements.
As a result of the surgery, the last component becomes knotted, and the
given manifold is its complement.
We found such surgery descriptions for 17 knot complements.

For the last 8 complements, we drilled out shortest geodesics in each
manifold and searched for any link $L$ in $S^3$ whose complement is
isometric.
The knot complement can then be expressed in terms of Dehn surgery on an
$(n-1)$-component sublink of $L$.
Again, Morwen Thistlethwaite provided tabulation data which proved
essential to find many of these links.
In several instances, we found SnapPea's bootstrapping manifolds, which
are link complements used to compute the Chern-Simons invariant.
To match surgery coefficients with the components of the $(n-1)$-component
sublink, we took advantage of a feature in Snap that keeps track of how
the surgery coefficients change when given several equivalent surgery
presentations.

Having obtained Dehn surgery descriptions for all 25 complements, we found
a sequence of Kirby calculus moves to reduce the surgery on $L$ to just
the single knotted component in $S^3$.



\section{Generalized twisted torus knots}
\vspace*{-0.5pt}
In this section we describe a generalization of twisted torus knots described in \cite{cdw99}. As motivation for generalizing twisted torus knots we illustrate a procedure to generate knots with high crossing number but whose complement has low volume. 

Let $M$ be a hyperbolic 3-manifold with n-cusps. Let $M(p_1,q_1)\ldots(p_n,q_n)$ denote the hyperbolic 3-manifold  obtained by $(p_i,q_i)$-filling the i-th cusp. Then it follows from work of Jorgensen and Thurston that 
$$ \mathrm{vol}(M(p_1,q_1)\ldots(p_n,q_n)) < \mathrm{vol}(M)$$

\begin{figure} 
\begin{center}
\psfrag{twists}{$n$ full left/right twists}
\psfrag{filling}{$\pm 1/n$}
\psfrag{volume}{3.663862= vol(Whitehead link) $>$ vol(knot complement) }
\includegraphics[width=5in]{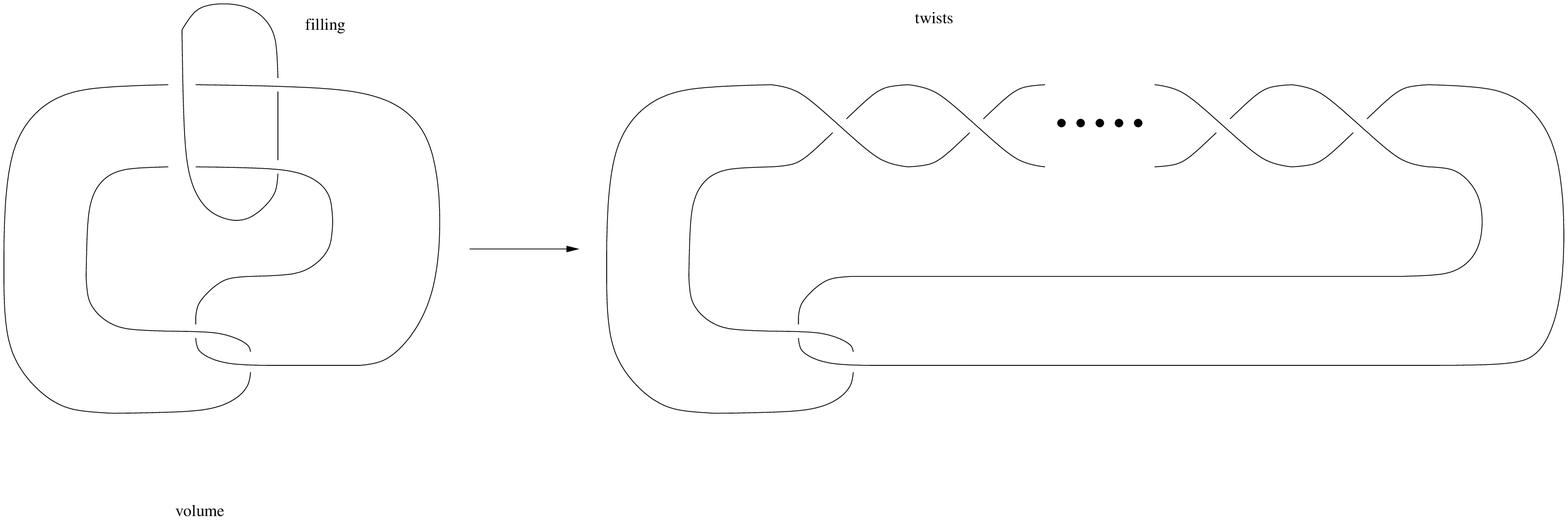}
\end{center}
\caption{$\pm 1/n$ surgery on Whitehead link} \label{whitehead}
\end{figure}

In particular, a knot obtained by $\pm 1/n$ filling on an unknotted component of a two component link will have high crossing number but volume bounded by the volume of the link complement, for example see Figure \ref{whitehead}. Similar ideas are used in \cite{lackenby} to obtain relations between a knot diagram and the volume. 

The twisted torus knot $T(p,q,r,s)$ can be seen as a $-1/s$ Dehn filling on the two component link with one component a $(p,q)$ torus knot, and the other component an unknot enclosing $r$ parallel strands. For an example, see Figure \ref{ttk}. Many twisted torus knots are hyperbolic and have interesting properties, especially small Seifert firbered surgeries \cite{cdw99}. We generalize twisted torus knots in two ways:

\begin{figure} 
\begin{center}
\psfrag{half}{$1/2$}
\includegraphics[height=3in]{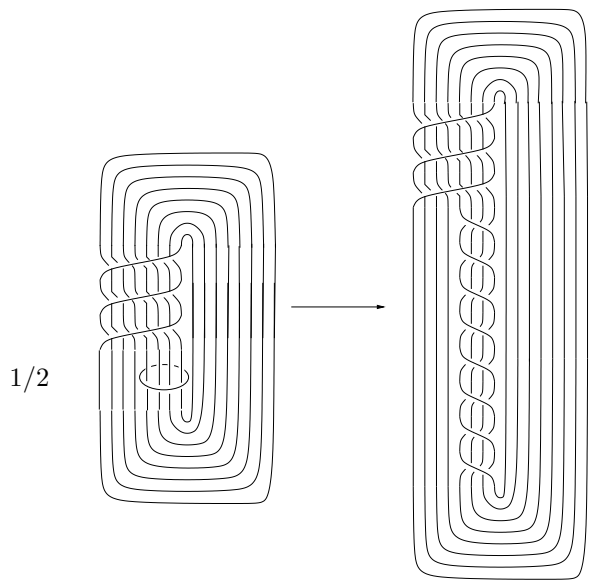}
\end{center}
\caption{Twisted torus knot T(8,3,4,-2)} \label{ttk}
\end{figure}

\begin{figure} 
\begin{center}
\psfrag{trefoil}{Trefoil}
\psfrag{markov}{Markov move}
\psfrag{twist2}{Two right twists on 3 strands}
\psfrag{t(5,3)}{$(5,2)$ torus knot}
\psfrag{3strands}{One right twist}
\psfrag{2strands}{2 left twists }
\psfrag{4strands}{One right twist}
\psfrag{move1}{First operation}
\psfrag{move2}{Second operation}
\includegraphics[height=3in]{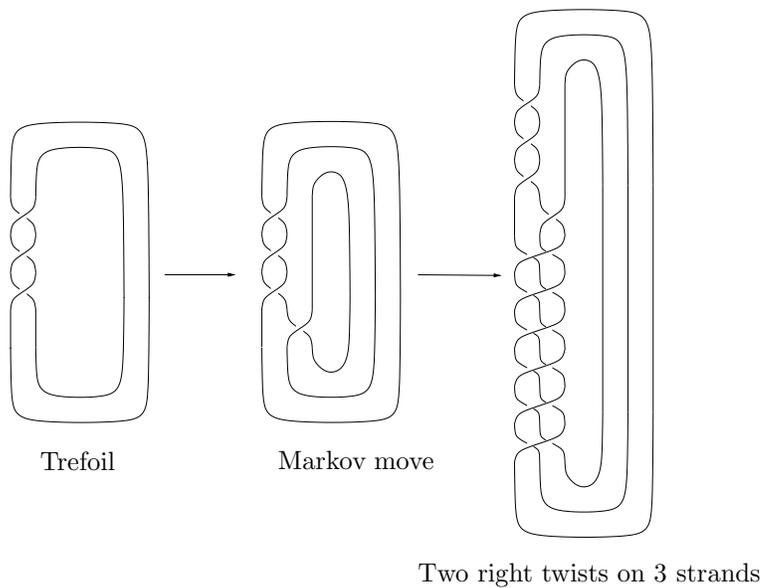}
\end{center}
\caption{Operation 1 for generalized twisted torus knots} \label{gttk1}
\end{figure}

\begin{figure} 
\begin{center}
\psfrag{trefoil}{Trefoil}
\psfrag{markov}{Markov move}
\psfrag{ttk}{Twisted torus knot}
\psfrag{T(5,3,3,1)}{ T$(5,2,3,1)$ }
\psfrag{3strands}{One right twist}
\psfrag{2strands}{Two left twists }
\psfrag{4strands}{One right twist}
\psfrag{move1}{First operation}
\psfrag{move2}{Second operation}
\includegraphics[height=3in]{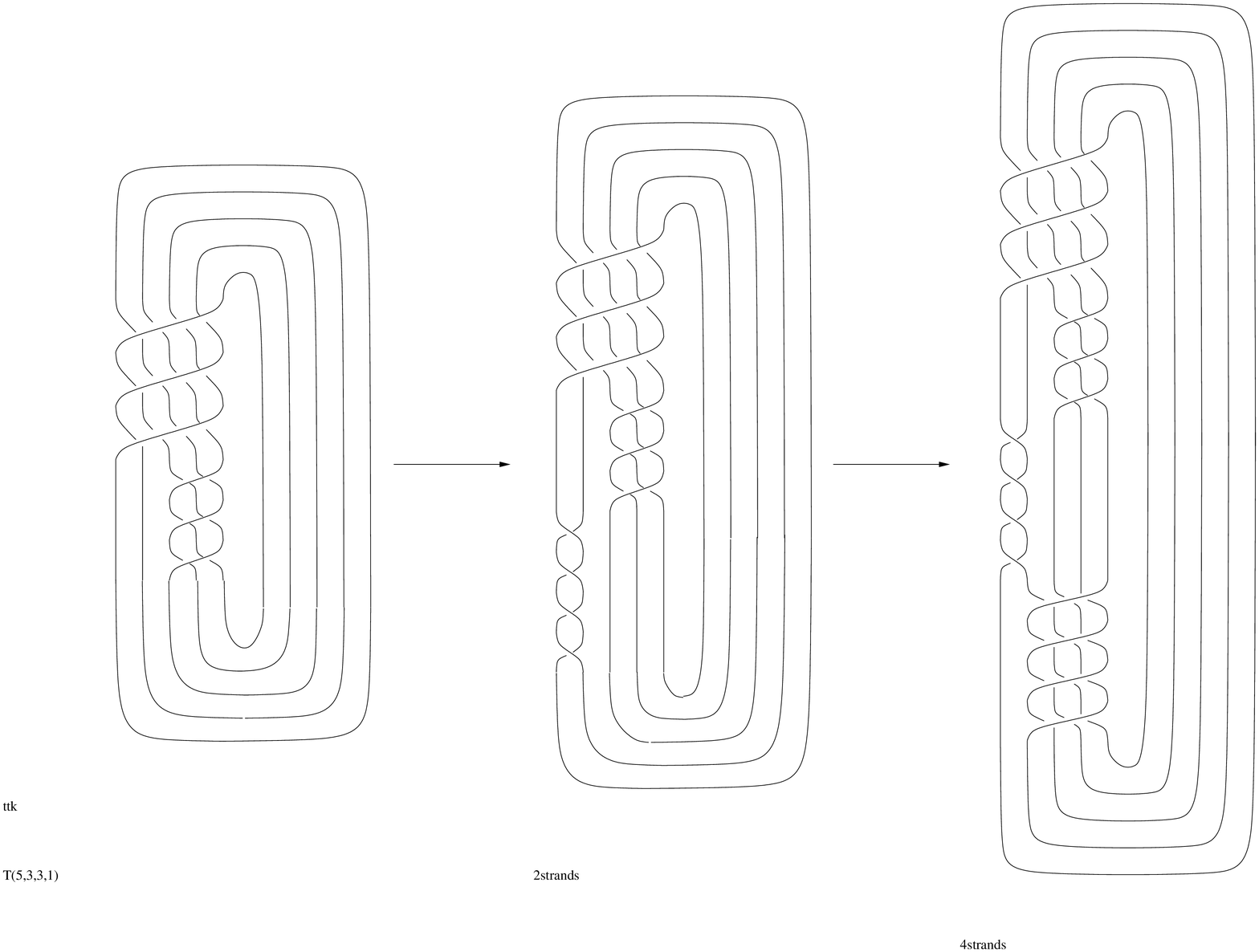}
\end{center}
\caption{Operation 2 for generalized twisted torus knots} \label{gttk2}
\end{figure}

\begin{enumerate}
\item Let $r > p$ by increasing the number of strands with Markov moves in the braid representation of the torus knot and then adding the $s$ twists on the $r$ strands. See Figure \ref{gttk1}.

\item Twist more than one group of strands. See Figure \ref{gttk2}.
\end{enumerate}

\begin{figure} 
\begin{center}
\psfrag{4left}{4 left}
\psfrag{3right}{3 right}
\psfrag{twists}{twists}
\psfrag{v2191}{knot corresponding to v2191}
\psfrag{v3234}{knot corresponding to v3234}
\includegraphics[height=3in]{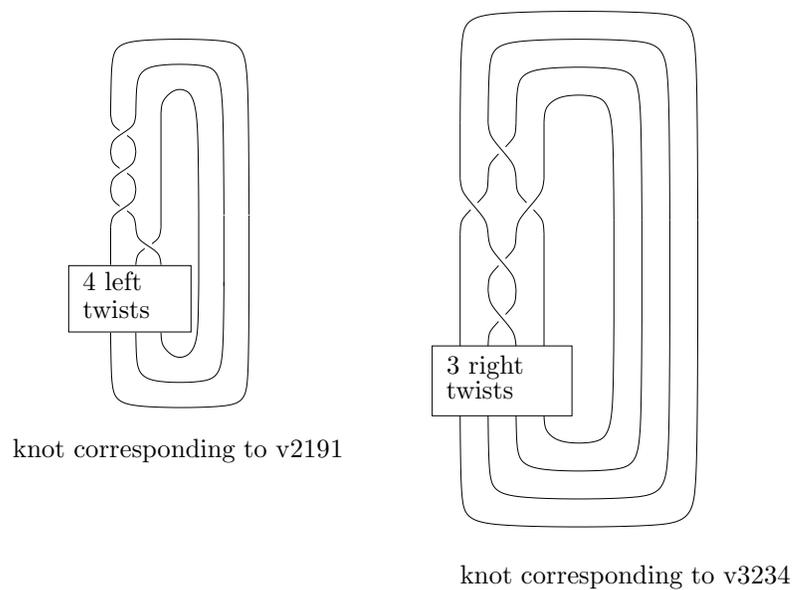}
\end{center}
\caption{Generalized twisted torus knots descriptions for some census knots} \label{gttkexamples}
\end{figure}

A knot obtained from a torus knot by a combination of the above moves is called a \textit{generalized twisted torus knot}. It would be interesting to study more properties of generalized twisted torus knots and ways in which they are similar to twisted torus knots. We have generalized twisted torus knot descriptions for some knots in our census, for example see Figure \ref{gttkexamples}.

\section{Kirby Calculus}
\vspace*{-0.5pt}
We use Kirby calculus on framed links to describe the remaining 25 knots which were not found by searching through Thistlethwaite's knot census or as twisted torus knots. The first step was to obtain a description of each knot complement as surgery on a framed link.

We follow \cite{saveliev99} and \cite{gs99} for Kirby Calculus. As an easy consequence of the Kirby moves let us describe some propositions which we use in our computations. 

\begin{figure} 
\begin{center}
\psfrag{n}{$n$}
\psfrag{nmp1}{$n \mp 1$}
\psfrag{pm1}{$\pm1$}
\psfrag{twist}{\scriptsize{left/right twist}}
\includegraphics[width=3in]{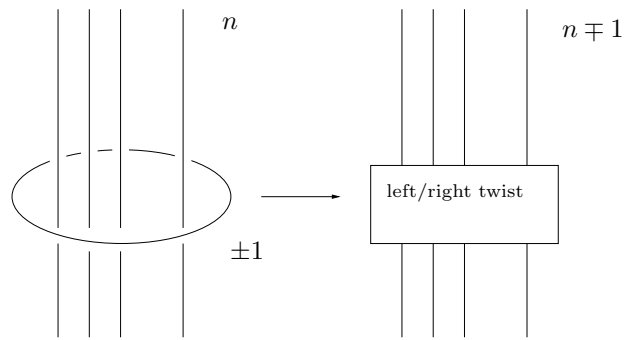}
\end{center}
\caption{Deleting component with framing $\pm 1$} \label{removepm1}
\end{figure}

{\bf Proposition~1.} {\it An unknot with framing $\pm 1$ can be deleted with the effect of giving all the arcs throught the unknot a full left/right twist and changing the framings by ading $\mp 1$ to each arc, assuming that they represent different components of $\mathcal{L}$. see Figure \ref{removepm1}. }


The case with linking number $2$ is described below.

\begin{figure} 
\begin{center}
\psfrag{n}{$n$}
\psfrag{nmp4}{$n\mp 4$}
\psfrag{n-4}{$n-4$}
\psfrag{pm1}{$\pm 1$}
\psfrag{twist}{left/right twist}
\includegraphics[width=3 in]{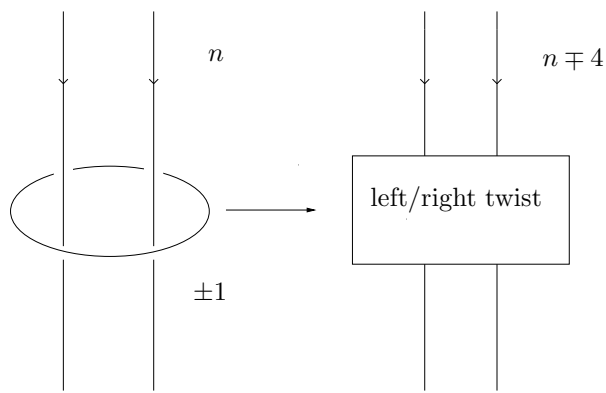}
\end{center}
\caption{Case with linking number 2} \label{linking2}
\end{figure}

{\bf Proposition~2.} {\it Assume a component of the link has linking number $2$ with an unknot with framing $\pm1$. Then the unknot can be deleted with the effect of giving the arcs a full left/right twist and changing the framings by adding $\mp 4$. See Figure \ref{linking2}. }


Another useful move which deals with rational framings is the \textit{slam-dunk}, whose dynamic name is due to T. Cochran.

\begin{figure} 
\begin{center}
\psfrag{z1}{$n \in \mathbb{Z}$}
\psfrag{q1}{$r \in \mathbb{Q}\cup \{ \infty \}$}
\psfrag{nr}{$\displaystyle{n-\frac{1}{r}}$}
\includegraphics[width=3in]{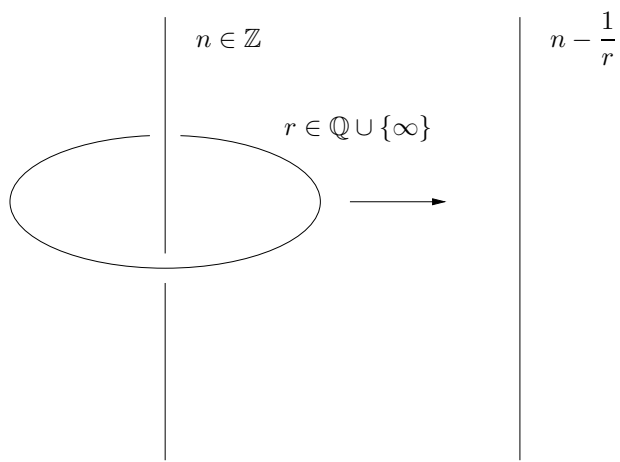}
\end{center}
\caption{Slam-dunk move} \label{slamdunk}
 \end{figure}

{\bf Proposition~3.} {\it  Suppose that one component $L_1$ is a meridian of another component $L_2$ and that the coefficients of these are $r \in \mathbb{Q}\cup \{ \infty \}$ and $n \in \mathbb{Z}$, respectively, then the unknot can be deleted with the effect of changing the framing of $L_2$ as $n-\frac{1}{r}$. See Figure \ref{slamdunk}. }


The operation shown in Figures 6-8 is called \textbf{blow-down}. The inverse operation is called \textbf{blow-up}.
\ \\


\begin{figure} 
\begin{center}
\psfrag{4}{$4$}
\psfrag{2/7}{$\displaystyle{\frac{2}{7}}$}
\psfrag{2}{$2$}
\psfrag{1}{$1$}
\psfrag{3}{$3$}
\psfrag{blowdown}{blow down three times}
\psfrag{slamdunk}{slam-dunk}
\includegraphics[width=5in]{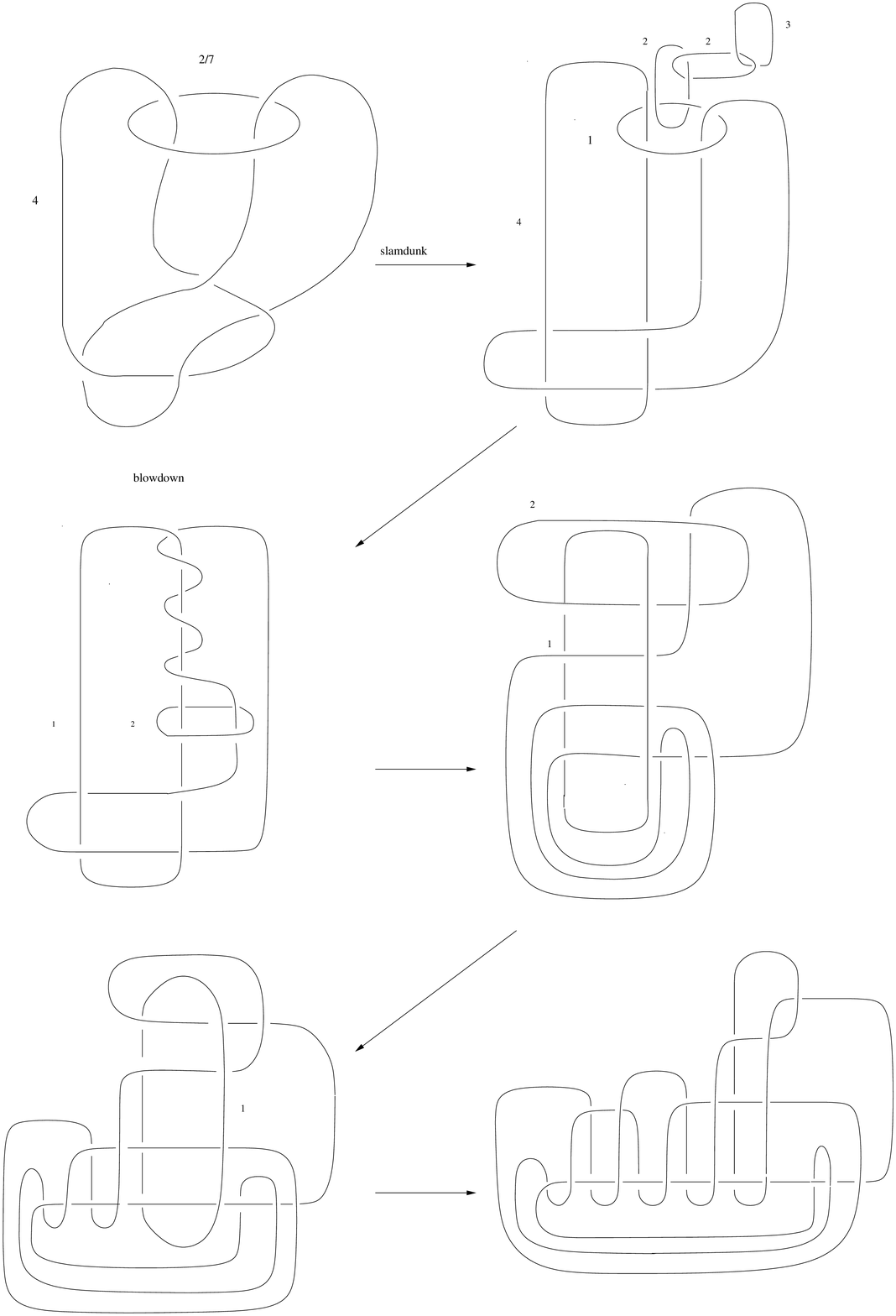}
\end{center}
\caption{Kirby calculus to obtain the knot for v656} \label{v656}
\end{figure}

{\bf Example.} In Figure \ref{v656}, we give an example of our technique for the knot corresponding the census manifold v656. The knot is obtained starting with a surgery description of v656 (the first picture in Figure \ref{v656}) and reducing components of the link using the propositions described above. The surgery description is obtained by drilling out shortest geodesics from v656 until the complement is isometric to a link complement with unknotted components. In this case, the link is $8_8^{3}$ in Rolfsen's table \cite{rolfsen}.
 
We expand $2/7$ as the following continued fraction:
$$\displaystyle{\frac{2}{7}=1-\frac{1}{2-\frac{1}{2-\frac{1}{3}}}}$$
Using the slam-dunk move above we can blow up to get integer surgery and carry out the computation as shown in Figure \ref{v656}.

\section{Tables}
 
We follow \cite{cdw99} for our notation and convention. Table 1 gives a list of knots whose complements can be decomposed into 7 ideal hyperbolic tetrahedra. We use the notation $\mathbf{k}7_m$ to indicate the $m^{th}$ knot in the list of knots made from $7$ tetrahedra. These knots are sorted in increasing volume. When there are multiple knots with the same volume we then sort by decreasing length of systole, the shortest closed geodesic in the complement. For example the knots $\mathbf{k}7_{95}$ and $\mathbf{k}7_{96}$ have the same volume but the knot $\mathbf{k}7_{95}$ has a longer systole and hence is listed earlier even though its number in SnapPea's census is higher. 

In Tables 2 and 3, we include the Jones polynomial for all knots,
including knots in \cite{cdw99}.  A glance at these polynomials reveals how
strikingly small they are compared with Jones polynomials of knots in
tables organized by crossing number.  Extreme examples are $\mathbf{k}7_{3},\
\mathbf{k}7_{4},\ \mathbf{k}7_{22},\ \mathbf{k}7_{23}$, whose polynomial spans are between 29 and 32,
but whose non-zero coefficients are all $\pm 1$ with absolute sum equal to 7.
The span of the Jones polynomial gives a lower bound for the crossing number, with equality if the knot is alternating.
The Jones polynomial with the largest coefficient average (about 2.64) is that
of $\mathbf{k}7_{103}$, a 10-crossing knot.  The Jones polynomial detects chirality,
so with help from Nathan Dunfield's python addition to SnapPea, we
have distinguished each knot from its mirror image according to the
orientation of the census manifold in SnapPea.  (In \cite{cdw99}, knots are
given only up to mirror image.)  The mirror image of T(p,q,r,s) is
denoted by T(p,-q,r,-s).  The Jones polynomials for $\mathbf{k}7_{61}$ and $\mathbf{k}7_{106}$
were computed by M. Ochiai using K2K version 1\_3\_1 \cite{ochiai}. The Jones polynomials for $\mathbf{k}7_{19}$ and $\mathbf{k}7_{124}$
were computed using knot diagrams simplified by Rob Scharein \cite{scharein}.  Computing
limitations left open the Jones polynomials for three knots.

In Table 4 we provide the Dowker-Thistlethwaite codes for the all 
the knots which do not have a twisted torus knot description. These are the 
 knots  described as ``See Below'' and the knots $C_k,\ Ca_k$ and $Cn_k$ in 
Table 1.

Here is a brief description of the table columns:
\begin{description}
\item[\textbf{OC: }] (old census) This gives the number of manifold as in \cite{hw89}. This nomenclature is also used in \textit{SnapPea}.
\item[\textbf{Volume: }]This is the hyperbolic volume of knot complement.
\item[\textbf{C-S: }]This is the Chern-Simons invariant of the knot complement. It is well defined for cusped hyperbolic 3-manifolds modulo 1/2. See \cite{meyerhoff92}.
\item[\textbf{Sym: }]This is the group of isometries of the knot complement. In the table \cyclic\ denotes the cyclic group of order 2 and \dihedral\ denotes the dihedral group of order 4.  Note that we do not get any other symmetry groups.\item[\textbf{SG: }]This is the length of the systole (shortest closed geodesic) in the knot complement. 
\item[\textbf{Description: }] Here we give a possible description of the knot as either a twisted torus knot or a knot in Rolfsen tables with 10 or fewer crossings or Thistlethwaite tables from 11 to 16 crossings. \\
An entry of the form $C_k$ for $C \leq 10$ indicates the $k^{th}$ knot of $C$ crossings in Rolfsen tables for knots with 10 or fewer crossings (\cite{rolfsen}). An entry of the form $Ca_k$ (respectively $Cn_k$) for $11 \leq C \leq 16$ indicates the $k^{th}$ alternating (respectively non-alternating) knot in Thistlethwaite's census of knots up to 16 crossings. An entry of the form $T(p,q,r,s)$ indicates a twisted torus knot as explained in section 2. The remaining knots were obtained by Kirby calculus on links with
unknotted components. Some of these have generalized twisted torus knot descriptions as explained in section 2.
\item[\textbf{Degree: }]The first integer gives the lowest degree and the second integer gives the highest degree of the Jones polynomial.
\item[\textbf{Jones polynomial: }]An entry of the form $(n,m)$ and $a_0+a_1+ \ldots +a_{n-m}$ corresponds to the polynomial $a_0t^n + a_1t^{n+1}+ \ldots + a_{n-m}t^m$.
\item[\textbf{C: }]This gives the number of crossings of the knot. When the knot is from Rolfsen's tables or Thistlethwaite's tables this is the minimal crossing number of the knot. For other knots this may not be the minimal crossing number.
\item[\textbf{DT code: }]This is the Dowker-Thistlethwaite code for the knot.
\end{description}

\newpage

 \newpage

 \begin{longtable}{|l||l|l|l|l|l|l|}

\caption{Knot complements with 7 tetrahedra}\\
  \hline
    & OC & Volume & C-S & Sym & SG & Description \\
  \hline 
  \endfirsthead

 \caption[]{\emph{continued}}\\
 \hline
   & OC & Volume & C-S & Sym & SG & Description \\
 \hline 
 \endhead

 \hline
 \endlastfoot

 $k7_{1}$ & 16 & 3.57388254 & 0.22001606 & \dihedral & 0.06 & $12a_{803}$ \\ 
\hline 
$k7_{2}$ & 25 & 3.58891392 & 0.03821170 & \dihedral & 0.05 & $13a_{3143}$ \\ 
\hline 
$k7_{3}$ & 82 & 3.62704008 & 0.11074361 & \cyclic & 0.02 & T(5,16,2,1) \\ 
\hline 
$k7_{4}$ & 114 & 3.63525119 & 0.17606301 & \cyclic & 0.02 & T(5,-17,3,-1) \\ 
\hline 
$k7_{5}$ & 165 & 3.90506656 & 0.08261094 & \cyclic & 0.10 & T(3,-17,2,-1) \\ 
\hline 
$k7_{6}$ & 220 & 4.02890933 & 0.12695546 & \cyclic & 0.02 & T(7,-17,2,1) \\ 
\hline 
$k7_{7}$ & 223 & 4.03440324 & 0.22601384 & \cyclic & 0.02 & T(7,18,2,1) \\ 
\hline 
$k7_{8}$ & 249 & 4.27985757 & 0.16931775 & \cyclic & 0.09 & T(6,-7,2,-1) \\ 
\hline 
$k7_{9}$ & 319 & 4.35277120 & 0.13945544 & \cyclic & 0.04 & T(6,-11,2,1) \\ 
\hline 
$k7_{10}$ & 321 & 4.35467010 & 0.01760473 & \cyclic & 0.04 & See below \\ 
\hline 
$k7_{11}$ & 329 & 4.35978253 & 0.19666614 & \cyclic & 0.04 & See below \\ 
\hline 
$k7_{12}$ & 330 & 4.36003045 & 0.21308713 & \cyclic & 0.13 & T(7,12,5,-1) \\ 
\hline 
$k7_{13}$ & 398 & 4.46284001 & 0.09756799 & \cyclic & 0.06 & T(7,-9,5,-1) \\ 
\hline 
$k7_{14}$ & 407 & 4.46826435 & 0.18469680 & \cyclic & 0.06 & T(8,-11,2,-1) \\ 
\hline 
$k7_{15}$ & 424 & 4.47895044 & 0.07613205 & \cyclic & 0.05 & T(13,3,11,-1) \\ 
\hline 
$k7_{16}$ & 434 & 4.48570291 & 0.14670089 & \cyclic & 0.04 & T(8,-13,2,1) \\ 
\hline 
$k7_{17}$ & 497 & 4.51180647 & 0.24439679 & \cyclic & 0.03 & T(12,19,11,-1) \\ 
\hline 
$k7_{18}$ & 521 & 4.53499561 & 0.21994497 & \cyclic & 0.14 & $15n_{80764}$ \\ 
\hline 
$k7_{19}$ & 535 & 4.55376662 & 0.19715482 & \cyclic & 0.13 & T(16,-5,9,1) \\ 
\hline 
$k7_{20}$ & 545 & 4.56838403 & 0.21271221 & \cyclic & 0.12 & T(5,11,4,-1) \\ 
\hline 
$k7_{21}$ & 554 & 4.57464504 & 0.22305545 & \cyclic & 0.21 & T(7,-9,2,-1) \\ 
\hline 
$k7_{22}$ & 570 & 4.59210710 & 0.20299839 & \cyclic & 0.08 & T(5,16,3,1) \\ 
\hline 
$k7_{23}$ & 573 & 4.59543707 & 0.08211521 & \cyclic & 0.08 & T(5,19,3,-1) \\ 
\hline 
$k7_{24}$ & 595 & 4.60805262 & 0.18359070 & \cyclic & 0.09 & See below \\ 
\hline 
$k7_{25}$ & 600 & 4.60992652 & 0.14728269 & \cyclic & 0.21 & $16n_{207543}$ \\ 
\hline 
$k7_{26}$ & 656 & 4.63782290 & 0.01834104 & \cyclic & 0.07 & See below \\ 
\hline 
$k7_{27}$ & 707 & 4.66145555 & 0.01027092 & \cyclic & 0.12 & T(4,17,2,-1) \\ 
\hline 
$k7_{28}$ & 709 & 4.66328880 & 0.22017834 & \cyclic & 0.20 & T(5,-9,4,-1) \\ 
\hline 
$k7_{29}$ & 715 & 4.66730782 & 0.05239840 & \cyclic & 0.05 & T(17,-5,14,1) \\ 
\hline 
$k7_{30}$ & 740 & 4.68848058 & 0.21080068 & \cyclic & 0.10 & T(4,19,2,1) \\ 
\hline 
$k7_{31}$ & 741 & 4.68947984 & 0.22307259 & \cyclic & 0.04 & T(18,-7,16,1) \\ 
\hline 
$k7_{32}$ & 759 & 4.70155763 & 0.21246644 & \cyclic & 0.12 & T(8,-11,2,1) \\ 
\hline 
$k7_{33}$ & 765 & 4.70790086 & 0.10488716 & \cyclic & 0.17 & T(7,12,2,-1) \\ 
\hline 
$k7_{34}$ & 830 & 4.75499484 & 0.06038348 & \cyclic & 0.18 & T(5,11,2,-1) \\ 
\hline 
$k7_{35}$ & 847 & 4.76458494 & 0.13833306 & \cyclic & 0.08 & T(7,-17,6,1) \\ 
\hline 
$k7_{36}$ & 912 & 4.80833192 & 0.21729973 & \cyclic & 0.15 & T(7,9,6,1) \\ 
\hline 
$k7_{37}$ & 939 & 4.84091784 & 0.16021452 & \cyclic & 0.16 & T(5,-14,2,-1) \\ 
\hline 
$k7_{38}$ & 945 & 4.84804861 & 0.11486758 & \cyclic & 0.13 & T(8,-13,2,-1) \\ 
\hline 
$k7_{39}$ & 959 & 4.85466334 & 0.13583195 & \cyclic & 0.30 & T(11,2,3,1) \\ 
\hline 
$k7_{40}$ & 960 & 4.85466334 & 0.05249861 & \cyclic & 0.30 & $12n_{235}$ \\ 
\hline 
$k7_{41}$ & 1063 & 4.93353007 & 0.05272917 & \cyclic & 0.13 & See below \\ 
\hline 
$k7_{42}$ & 1077 & 4.94762401 & 0.19061542 & \cyclic & 0.12 & T(7,17,2,1) \\ 
\hline 
$k7_{43}$ & 1109 & 4.97329565 & 0.20882427 & \cyclic & 0.12 & T(7,18,2,-1) \\ 
\hline 
$k7_{44}$ & 1126 & 4.99345656 & 0.11840638 & \cyclic & 0.09 & See below \\ 
\hline 
$k7_{45}$ & 1217 & 5.11484146 & 0.16966684 & \dihedral & 0.14 & $10_2$ \\ 
\hline 
$k7_{46}$ & 1243 & 5.14020656 & 0.07876069 & \dihedral & 0.12 & $11a_{364}$ \\ 
\hline 
$k7_{47}$ & 1269 & 5.16583510 & 0.11994300 & \cyclic & 0.17 & T(6,-7,4,-1) \\ 
\hline 
$k7_{48}$ & 1300 & 5.18904995 & 0.14371112 & \cyclic & 0.27 & T(7,10,3,-1) \\ 
\hline 
$k7_{49}$ & 1348 & 5.21927466 & 0.12925600 & \cyclic & 0.29 & T(4,5,3,-3) \\ 
\hline 
$k7_{50}$ & 1359 & 5.22537658 & 0.16444233 & \cyclic & 0.11 & T(8,-5,7,2) \\ 
\hline 
$k7_{51}$ & 1392 & 5.24108057 & 0.00759382 & \cyclic & 0.26 & See below \\ 
\hline 
$k7_{52}$ & 1423 & 5.25717856 & 0.22317376 & \cyclic & 0.16 & T(6,-11,4,1) \\ 
\hline 
$k7_{53}$ & 1425 & 5.25881410 & 0.17118622 & \cyclic & 0.37 & T(10,-3,4,-1) \\ 
\hline 
$k7_{54}$ & 1434 & 5.26280499 & 0.03696783 & \cyclic & 0.26 & See below \\ 
\hline 
$k7_{55}$ & 1547 & 5.33947334 & 0.01495549 & \cyclic & 0.25 & See below \\ 
\hline 
$k7_{56}$ & 1565 & 5.34928267 & 0.11217244 & \cyclic & 0.31 & T(4,17,3,-1) \\ 
\hline 
$k7_{57}$ & 1620 & 5.37586245 & 0.12860526 & \cyclic & 0.29 & T(4,15,3,1) \\ 
\hline 
$k7_{58}$ & 1628 & 5.37781448 & 0.17721031 & \cyclic & 0.15 & T(7,-10,4,-1) \\ 
\hline 
$k7_{59}$ & 1690 & 5.40911134 & 0.18400493 & \cyclic & 0.16 & See below \\ 
\hline 
$k7_{60}$ & 1709 & 5.41958732 & 0.18626556 & \cyclic & 0.10 & T(8,-5,7,-1) \\ 
\hline 
$k7_{61}$ & 1716 & 5.42554802 & 0.15249052 & \cyclic & 0.19 & T(14,5,12,-1) \\ 
\hline 
$k7_{62}$ & 1718 & 5.42639797 & 0.17745037 & \cyclic & 0.11 & T(10,3,8,1) \\ 
\hline 
$k7_{63}$ & 1728 & 5.43451894 & 0.14955510 & \cyclic & 0.24 & T(7,5,6,-3) \\ 
\hline 
$k7_{64}$ & 1810 & 5.47368625 & 0.12109788 & \cyclic & 0.27 & T(5,17,4,1) \\ 
\hline 
$k7_{65}$ & 1811 & 5.47375747 & 0.21730448 & \cyclic & 0.24 & See below \\ 
\hline 
$k7_{66}$ & 1832 & 5.48650517 & 0.05905336 & \cyclic & 0.13 & T(7,4,6,-3) \\ 
\hline 
$k7_{67}$ & 1839 & 5.48907036 & 0.04062381 & \cyclic & 0.13 & See below \\ 
\hline 
$k7_{68}$ & 1915 & 5.52030806 & 0.05903973 & \cyclic & 0.29 & T(3,-11,2,-2) \\ 
\hline 
$k7_{69}$ & 1921 & 5.52226249 & 0.02153669 & \cyclic & 0.11 & T(10,3,8,-2) \\ 
\hline 
$k7_{70}$ & 1935 & 5.52559924 & 0.05425344 & \cyclic & 0.69 & $13n_{192}$ \\ 
\hline 
$k7_{71}$ & 1940 & 5.52812617 & 0.17689247 & \cyclic & 0.13 & T(7,-10,6,-1) \\ 
\hline 
$k7_{72}$ & 1964 & 5.54249469 & 0.11539683 & \cyclic & 0.19 & $16n_{679988}$ \\ 
\hline 
$k7_{73}$ & 1966 & 5.54388606 & 0.16468534 & \cyclic & 0.09 & T(21,-8,7,1) \\ 
\hline 
$k7_{74}$ & 1971 & 5.54788269 & 0.13580814 & \cyclic & 0.24 & See below \\ 
\hline 
$k7_{75}$ & 1980 & 5.55020317 & 0.23223326 & \cyclic & 0.54 & T(5,-8,3,-1) \\ 
\hline 
$k7_{76}$ & 1986 & 5.55193476 & 0.17160973 & \cyclic & 0.12 & T(7,4,6,2) \\ 
\hline 
$k7_{77}$ & 2001 & 5.55993266 & 0.24388215 & \dihedral & 0.20 & See below \\ 
\hline 
$k7_{78}$ & 2024 & 5.56671464 & 0.23687581 & \cyclic & 0.21 & T(7,5,6,2) \\ 
\hline 
$k7_{79}$ & 2090 & 5.59955942 & 0.19922862 & \cyclic & 0.13 & T(21,-8,17,1) \\ 
\hline 
$k7_{80}$ & 2166 & 5.63877295 & 0.19098985 & \cyclic & 0.59 & $10_{161}$ \\ 
\hline 
$k7_{81}$ & 2191 & 5.64969858 & 0.18611119 & \cyclic & 0.21 & T(3,13,2,-2) \\ 
\hline 
$k7_{82}$ & 2215 & 5.66111586 & 0.24338756 & \cyclic & 0.25 & T(5,-13,4,-1) \\ 
\hline 
$k7_{83}$ & 2217 & 5.66193313 & 0.13327962 & \cyclic & 0.54 & T(11,2,7,-1) \\ 
\hline 
$k7_{84}$ & 2257 & 5.68694120 & 0.06713052 & \dihedral & 0.12 & See below \\ 
\hline 
$k7_{85}$ & 2272 & 5.69302109 & 0.08333333 & \cyclic & 0.55 & $13n_{469}$ \\ 
\hline 
$k7_{86}$ & 2284 & 5.69844175 & 0.10029838 & \dihedral & 0.29 & $9_5$ \\ 
\hline 
$k7_{87}$ & 2290 & 5.69905559 & 0.09594144 & \cyclic & 0.35 & T(5,-9,3,-1) \\ 
\hline 
$k7_{88}$ & 2325 & 5.71567170 & 0.03475951 & \cyclic & 0.24 & T(7,3,5,3) \\ 
\hline 
$k7_{89}$ & 2362 & 5.73210479 & 0.07217067 & \dihedral & 0.27 & $10_3$ \\ 
\hline 
$k7_{90}$ & 2384 & 5.75078061 & 0.06729087 & \cyclic & 0.50 & T(3,-7,2,-4) \\ 
\hline 
$k7_{91}$ & 2488 & 5.81712969 & 0.12620932 & \dihedral & 0.21 & $10_4$ \\ 
\hline 
$k7_{92}$ & 2508 & 5.82963473 & 0.09558485 & \cyclic & 0.47 & $14n_{14254}$ \\ 
\hline 
$k7_{93}$ & 2520 & 5.83875463 & 0.15594177 & \dihedral & 0.19 & $11a_{342}$ \\ 
\hline 
$k7_{94}$ & 2543 & 5.85424272 & 0.09614515 & \cyclic & 0.86 & $12n_{582}$ \\ 
\hline 
$k7_{95}$ & 2553 & 5.86053930 & 0.22027574 & \cyclic & 0.62 & $10_{128}$ \\ 
\hline 
$k7_{96}$ & 2552 & 5.86053930 & 0.13694240 & \cyclic & 0.28 & $11n_{57}$ \\ 
\hline 
$k7_{97}$ & 2623 & 5.90408586 & 0.21800323 & \cyclic & 0.74 & $9_{43}$ \\ 
\hline 
$k7_{98}$ & 2624 & 5.90408586 & 0.19866344 & \cyclic & 0.25 & $12n_{243}$ \\ 
\hline 
$k7_{99}$ & 2642 & 5.91674574 & 0.19843862 & \dihedral & 1.11 & $12n_{725}$ \\ 
\hline 
$k7_{100}$ & 2743 & 5.98392087 & 0.18936943 & \dihedral & 0.36 & See below \\ 
\hline 
$k7_{101}$ & 2759 & 5.99627882 & 0.11017406 & \cyclic & 0.42 & T(7,11,3,-1) \\ 
\hline 
$k7_{102}$ & 2851 & 6.07547284 & 0.23810941 & \cyclic & 0.46 & T(4,5,2,2) \\ 
\hline 
$k7_{103}$ & 2858 & 6.08323484 & 0.10630991 & \dihedral & 0.30 & $10_8$ \\ 
\hline 
$k7_{104}$ & 2861 & 6.08444991 & 0.15830176 & \cyclic & 0.74 & $13n_{3523}$ \\ 
\hline 
$k7_{105}$ & 2869 & 6.08730233 & 0.00297891 & \cyclic & 0.35 & See below \\ 
\hline 
$k7_{106}$ & 2871 & 6.08754805 & 0.14502692 & \cyclic & 0.96 & T(14,3,8,-1) \\ 
\hline 
$k7_{107}$ & 2888 & 6.09832261 & 0.05473501 & \cyclic & 0.36 & See below \\ 
\hline 
$k7_{108}$ & 2894 & 6.10280316 & 0.22005821 & \dihedral & 0.29 & $11a_{358}$ \\ 
\hline 
$k7_{109}$ & 2900 & 6.10811906 & 0.24457111 & \cyclic & 0.28 & See below \\ 
\hline 
$k7_{110}$ & 2925 & 6.12496920 & 0.09213713 & \cyclic & 0.76 & See below \\ 
\hline 
$k7_{111}$ & 2930 & 6.12852645 & 0.23594903 & \cyclic & 0.47 & T(7,-10,3,-1) \\ 
\hline 
$k7_{112}$ & 3070 & 6.23517885 & 0.00706941 & \cyclic & 0.34 & T(8,-5,6,3) \\ 
\hline 
$k7_{113}$ & 3081 & 6.24376275 & 0.09862176 & \cyclic & 0.35 & T(4,7,2,-2) \\ 
\hline 
$k7_{114}$ & 3093 & 6.25865905 & 0.21083423 & \cyclic & 0.43 & $16n_{245346}$ \\ 
\hline 
$k7_{115}$ & 3105 & 6.26514111 & 0.03296114 & \cyclic & 0.51 & See below \\ 
\hline 
$k7_{116}$ & 3169 & 6.31703999 & 0.21555521 & \cyclic & 0.42 & See below \\ 
\hline 
$k7_{117}$ & 3195 & 6.34376089 & 0.20313717 & \cyclic & 0.74 & $12n_{309}$ \\ 
\hline 
$k7_{118}$ & 3199 & 6.34557931 & 0.17015657 & \dihedral & 0.72 & See below \\ 
\hline 
$k7_{119}$ & 3234 & 6.36297460 & 0.22039529 & \cyclic & 0.57 & See below \\ 
\hline 
$k7_{120}$ & 3310 & 6.44353738 & 0.12055587 & \dihedral & 0.84 & $7_5$ \\ 
\hline 
$k7_{121}$ & 3320 & 6.45137056 & 0.11580073 & \cyclic & 1.11 & $12n_{475}$ \\ 
\hline 
$k7_{122}$ & 3335 & 6.47226750 & 0.01184030 & \cyclic & 0.62 & T(15,4,10,-1) \\ 
\hline 
$k7_{123}$ & 3354 & 6.51377663 & 0.13646474 & \cyclic & 0.29 & T(8,5,6,2) \\ 
\hline 
$k7_{124}$ & 3356 & 6.52092765 & 0.09834517 & \dihedral & 0.65 & T(15,4,6,-1) \\ 
\hline 
$k7_{125}$ & 3423 & 6.59122061 & 0.18267803 & \cyclic & 0.56 & $15n_{103488}$ \\ 
\hline 
$k7_{126}$ & 3482 & 6.72238568 & 0.22973136 & \cyclic & 0.61 & T(8,-3,4,-2) \\ 
\hline 
$k7_{127}$ & 3505 & 6.78371352 & 0.24180460 & \dihedral & 1.34 & $8_{21}$ \\ 
\hline 
$k7_{128}$ & 3536 & 6.90911013 & 0.03063918 & \cyclic & 0.86 & $11n_{49}$ \\ 
\hline 
$k7_{129}$ & 3541 & 6.92263444 & 0.03204283 & \cyclic & 0.46 & See below \\ 
\hline

 \end{longtable}
 \newpage

 \begin{longtable}{|l||l|p{3.55in}|}

 \caption{Jones polynomial of Knots with $\leq 6$  tetrahedra}\\
  \hline
   &Degree& Jones polynomial \\
  \hline 
  \endfirsthead

 \caption[]{\emph{continued}}\\
 \hline
  &Degree& Jones polynomial \\
 \hline 
 \endhead

 \hline
 \endlastfoot

 $\mathbf{k}2_{1} $ &(-2, 2)& $1-1+1-1+1  $\\ 
\hline 
$\mathbf{k}3_{1} $ &(5, 13)& $1+0+1+0+0+0-1+1-1  $\\ 
\hline 
$\mathbf{k}3_{2} $ &(-6, -1)& $-1+1-1+2-1+1  $\\ 
\hline 
$\mathbf{k}4_{1} $ &(-4, 2)& $1-1+1-2+2-1+1  $\\ 
\hline 
$\mathbf{k}4_{2} $ &(-8, -1)& $-1+1-1+2-2+2-1+1  $\\ 
\hline 
$\mathbf{k}4_{3} $ &(-19, -8)& $-1+1-1+0+0+0+0+0+0+1+0+1  $\\ 
\hline 
$\mathbf{k}4_{4} $ &(-23, -10)& $1-2+1-1+0+0-1+1-1+1+0+1+0+1  $\\
\hline 
$\mathbf{k}5_{1} $ &(-24, -11)& $1-1+0-1+0-1+0+0+0+1+0+1+0+1  $\\ 
\hline 
$\mathbf{k}5_{2} $ &(-6, 2)& $1-1+1-2+2-2+2-1+1  $\\ 
\hline 
$\mathbf{k}5_{3} $ &(-10, -1)& $-1+1-1+2-2+2-2+2-1+1  $\\ 
\hline 
$\mathbf{k}5_{4} $ &(15, 31)& $1+0+1+0+1+0+0+0+0-1+0-1+0+0+0+1-1  $\\ 
\hline 
$\mathbf{k}5_{5} $ &(11, 25)& $1+0+1+0+0+0+0+0+0+0+0+0-1+1-1  $\\ 
\hline 
$\mathbf{k}5_{6} $ &(13, 27)& $1+0+1+0+1+0+0+0+0-1+0-1+0-1+1  $\\ 
\hline 
$\mathbf{k}5_{7} $ &(-28, -13)& $-1+1+0+0+0-1+0-1+0+0+0+1+0+1+0+1  $\\ 
\hline 
$\mathbf{k}5_{8} $ &(-3, 3)& $1-1+1-1+1-1+1  $\\ 
\hline 
$\mathbf{k}5_{9} $ &(-7, -2)& $-1+1-1+1+0+1  $\\ 
\hline 
$\mathbf{k}5_{10} $ &(9, 20)& $1+0+1+0+1-1+1-2+1-1+1-1  $\\ 
\hline 
$\mathbf{k}5_{11} $ &(6, 16)& $1+0+1+0+0+0-1+1-1+1-1  $\\ 
\hline 
$\mathbf{k}5_{12} $ &(-1, 5)& $-1+2-1+2-1+1-1  $\\ 
\hline 
$\mathbf{k}5_{13} $ &(-2, 5)& $1+0+0+0-1+1-1+1  $\\ 
\hline 
$\mathbf{k}5_{14} $ &(11, 25)& $1+0+1+0+1-1+1-1+0+0-1+1-2+2-1  $\\ 
\hline 
$\mathbf{k}5_{15} $ &(12, 24)& $1+0+1+0+1-1+1-1+0+0+0+0-1  $\\ 
\hline 
$\mathbf{k}5_{16} $ &(16, 33)& $1+0+1+0+1+0+0+0+0-1+0+0-1+1-1+1-2+1  $\\ 
\hline 
$\mathbf{k}5_{17} $ &(18, 37)& $1+0+1+0+1+0+0+0+0+0-1+1-2+1-1+0+0-1+2-1  $\\ 
\hline 
$\mathbf{k}5_{18} $ &(16, 35)& $1+0+1+0+1-1+1-1+1-1+1-1+0+0-1+1-2+2-2+1  $\\ 
\hline 
$\mathbf{k}5_{19} $ &(-1, 5)& $1-1+2-2+2-2+1  $\\ 
\hline 
$\mathbf{k}5_{20} $ &(2, 9)& $1-1+2-2+3-2+1-1  $\\ 
\hline 
$\mathbf{k}5_{21} $ &(0, 6)& $2-1+1-2+1-1+1  $\\ 
\hline 
$\mathbf{k}5_{22} $ &(4, 12)& $1+0+1+0-1+1-1+1-1  $\\ 
\hline 
$\mathbf{k}6_{1} $ &(-2, 8)& $1-1+2-2+2-2+2-2+1-1+1  $\\ 
\hline 
$\mathbf{k}6_{2} $ &(-12, -1)& $-1+1-1+2-2+2-2+2-2+2-1+1  $\\ 
\hline 
$\mathbf{k}6_{3} $ &(21, 42)& $1+0+1+0+1+0+0+0+0+0+0+0+0-1+0-1+0+0+0+0-1+1  $\\ 
\hline 
$\mathbf{k}6_{4} $ &(25, 49)& $1+0+1+0+1+0+0+0+0+0+0+0+0+0-1+0-1+0+0+0+0+0+0+1-1  $\\ 
\hline 
$\mathbf{k}6_{5} $ &(14, 31)& $1+0+1+0+0+0+0+0+0+0+0+0+0+0+0-1+1-1  $\\ 
\hline 
$\mathbf{k}6_{6} $ &(26, 50)& $1+0+1+0+1+0+1+0+0+0+0-1+0-1+0-1+0+0-1+1-1+1+0+1-1  $\\ 
\hline 
$\mathbf{k}6_{7} $ &(-59, -31)& $1-1+0+0-1+1-1+1+0+0+0+0-1+0-1+0-1+0+0+0+0+0+1+0+1+0+1+0+1  $\\ 
\hline 
$\mathbf{k}6_{8} $ &(-10, -2)& $-1+1-1+1+0+0+0+0+1  $\\ 
\hline 
$\mathbf{k}6_{9} $ &(-8, 2)& $1-1+1-1+0+0+0+0+0+0+1  $\\ 
\hline 
$\mathbf{k}6_{10} $ &(1, 9)& $1-1+1+0+0+1-1+1-1  $\\ 
\hline 
$\mathbf{k}6_{11} $ &(-35, -17)& $1-1+0+0-1+1-2+1-1+0+0+0+0+0+1+0+1+0+1  $\\ 
\hline 
$\mathbf{k}6_{12} $ &(23, 45)& $1+0+1+0+1+0+0+0+0+0+0+0+0+0-1+0-1+0+0+0+0-1+1  $\\ 
\hline 
$\mathbf{k}6_{13} $ &(23, 46)& $1+0+1+0+1+0+0+0+0+0+0+0+0-1+0-1+0+0+0+0+0+0+1-1  $\\ 
\hline 
$\mathbf{k}6_{14} $ &(17, 37)& $1+0+1+0+1-1+1-1+1-1+1-1+0+0-1+1-2+2-2+2-1  $\\ 
\hline 
$\mathbf{k}6_{15} $ &(24, 48)& $1+0+1+0+1+0+1-1+1-1+1-2+1-1+0-1+1-1+0+0-1+1-1+2-1  $\\ 
\hline 
$\mathbf{k}6_{16} $ &(28, 54)& $1+0+1+0+1+0+1+0+0+0+0-1+0-1+0-1+0+0+0+0+1-1+1-1+0-1+1  $\\ 
\hline 
$\mathbf{k}6_{17} $ &(22, 47)& $1+0+1+0+1-1+1-1+1-1+1-1+1-1+1-1+0+0-1+1-2+2-2+2-2+1  $\\ 
\hline 
$\mathbf{k}6_{18} $ &(27, 53)& $1+0+1+0+1+0+1-1+1-1+1-1+0+0-1+0-1+1-2+2-2+1+0+0+1-2+1  $\\ 
\hline 
$\mathbf{k}6_{19} $ &(7, 19)& $1+0+1+0+0+0-1+1-1+1-1+1-1  $\\ 
\hline 
$\mathbf{k}6_{20} $ &(-4, 4)& $-1+1-1+2-1+2-1+1-1  $\\ 
\hline 
$\mathbf{k}6_{21} $ &(-55, -29)& $-1+1+0+0+1-1+1-1+0+0-1+0-1+0-1+0+0+0+0+0+1+0+1+0+1+0+1  $\\ 
\hline 
$\mathbf{k}6_{22} $ &(-2, 2)& $1-1+1-1+1  $\\ 
\hline 
$\mathbf{k}6_{23} $ &(-8, 0)& $1-2+2-3+3-2+2-1+1  $\\ 
\hline 
$\mathbf{k}6_{24} $ &(-12, -3)& $-1+1-2+3-3+3-2+2-1+1  $\\ 
\hline 
$\mathbf{k}6_{25} $ &(15, 30)& $1+0+1+0+1-1+1-1+1-1+1-2+1-1+1-1  $\\ 
\hline 
$\mathbf{k}6_{26} $ &(-10, -2)& $-1+1-1+1+0+0+0+0+1  $\\ 
\hline 
$\mathbf{k}6_{27} $ &(-34, -18)& $-1+0+0+0+0-1+1-1+1-1+1-1+1+0+1+0+1  $\\ 
\hline 
$\mathbf{k}6_{28} $ &(1, 8)& $1-2+3-2+3-2+1-1  $\\ 
\hline 
$\mathbf{k}6_{29} $ &(9, 22)& $1+0+1+0+0+0+0+0+0-1+1-1+1-1  $\\ 
\hline 
$\mathbf{k}6_{30} $ &(-55, -28)& $-1+2-1+0+0+0+0+0-1+1-2+1-1+0+0+0+0+0+0+0+0+0+0+ 1+0+1+0+1  $\\ 
\hline 
$\mathbf{k}6_{31} $ &(-22, -9)& $1-1+0-1+0+0-1+1-1+1+0+1+0+1  $\\ 
\hline 
$\mathbf{k}6_{32} $ &(-41, -20)& $1-2+1+0+0+0-1+0+0-1+1-1+0+0+0+0+0+1+0+1+0+1  $\\ 
\hline 
$\mathbf{k}6_{33} $ &(-7, 0)& $-1+1-1+2-1+1-1+1  $\\ 
\hline 
$\mathbf{k}6_{34} $ &(-4, 4)& $1-1+2-3+3-3+2-1+1  $\\ 
\hline 
$\mathbf{k}6_{35} $ &(-57, -29)& $-1+2-1+0+0+0+1-1+0+0-1+0+0-1+0+0+0-1+1-1+1-1+1+0+1+0+1+0+1  $\\ 
\hline 
$\mathbf{k}6_{36} $ &(26, 51)& $1+0+1+0+1+0+0+0+0+0+0+0+0+0-1+0+0-1+1-1+0+0+0+1-2+1  $\\ 
\hline 
$\mathbf{k}6_{37} $ &(-4, 5)& $-1+1-1+1+1+0+1-1+1-1  $\\ 
\hline 
$\mathbf{k}6_{38} $ &(7, 18)& $1+0+1+0+0+0+0-1+1-1+1-1  $\\ 
\hline 
$\mathbf{k}6_{39} $ &(13, 29)& $1+0+1+0+1-1+1-1+0+0+0+0+0-1+1-2+1  $\\ 
\hline 
$\mathbf{k}6_{40} $ &(-3, 5)& $1-1+2-3+3-3+3-2+1  $\\ 
\hline 
$\mathbf{k}6_{41} $ &(-11, -2)& $-1+1-2+3-3+4-3+2-1+1  $\\ 
\hline 
$\mathbf{k}6_{42} $ &(0, 8)& $1+0+1-1+1-2+1-1+1  $\\ 
\hline 
$\mathbf{k}6_{43} $ &(-3, 3)& $-1+2-2+3-2+2-1  $\\ 
\hline

 \end{longtable}

 \newpage

 \begin{longtable}{|l||l|p{3.55in}|}

 \caption{Jones polynomial of Knots with 7 tetrahedra}\\
  \hline
   &Degree& Jones polynomial \\
  \hline 
  \endfirsthead

 \caption[]{\emph{continued}}\\
 \hline
  &Degree& Jones polynomial \\
 \hline 
 \endhead

 \hline
 \endlastfoot

 $\mathbf{k}7_{1} $ &(-10, 2)& $1-1+1-2+2-2+2-2+2-2+2-1+1  $\\ 
\hline 
$\mathbf{k}7_{2} $ &(-14, -1)& $-1+1-1+2-2+2-2+2-2+2-2+2-1+1  $\\ 
\hline 
$\mathbf{k}7_{3} $ &(31, 60)& $1+0+1+0+1+0+0+0+0+0+0+0+0+0+0+0+0+0-1+0-1+0+0+0+0+0+0+0-1+1  $\\ 
\hline 
$\mathbf{k}7_{4} $ &(-67, -35)& $-1+1+0+0+0+0+0+0+0+0+0-1+0-1+0+0+0+0+0+0+0+0+0+ 0+0+0+0+0+1+0+1+0+1  $\\ 
\hline 
$\mathbf{k}7_{5} $ &(-37, -17)& $-1+1-1+0+0+0+0+0+0+0+0+0+0+0+0+0+0-0+1+0+1  $\\ 
\hline 
$\mathbf{k}7_{6} $ &(-86, -47)& $-1+1+0+0+0+0+1-1+1-1+0+0+0+0+0+0+0-1+0-1+0-1+0+0+0+0+0+0+0+0+0+0+0+1+0+1+0+1+0+1  $\\ 
\hline 
$\mathbf{k}7_{7} $ &(52, 95)& $1+0+1+0+1+0+1+0+0+0+0+0+0+0+0+0+0-0+0-1+0-1+0-1+0+0+0+0+0+0+0+0+0+1-1+1-1+0+0+0+0+0-1+1  $\\ 
\hline 
$\mathbf{k}7_{8} $ &(-32, -16)& $1-1+1-2+1-2+1-2+1-1+1+0+1+0+1+0+1  $\\ 
\hline 
$\mathbf{k}7_{9} $ &(-45, -24)& $-1+1-1+2-2+2-2+1-2+1-2+1-1+1-1+1+0+1+0+1+0+1  $\\ 
\hline 
$\mathbf{k}7_{10} $ &(2, 13)& $1+0+0+0+0+0+0+0+1-1+1-1  $\\ 
\hline 
$\mathbf{k}7_{11} $ &(-2, 11)& $1+0+0+0+0+0+0+0+0+0-1+1-1+1  $\\ 
\hline 
$\mathbf{k}7_{12} $ &(23, 45)& $1+0+1+0+1+0+1+0+0+0-1+0-1+0-1+0+0+0+0+0+0-1+1  $\\ 
\hline 
$\mathbf{k}7_{13} $ &(-64, -34)& $-1+1+0+0+0+0+0+0+0+0+0+0+0-1+0-1+0-1+0+0+0+0+0+0+1+0+1+0+1+0+1  $\\ 
\hline 
$\mathbf{k}7_{14} $ &(-68, -36)& $1-2+1-1+1+0+0+0+0+0+0+0+0+0+0-1+0-1+0-1+0-1+1-1+1+0+1+0+1+0+1+0+1  $\\ 
\hline 
$\mathbf{k}7_{15} $ &(-63, -33)& $-1+2-2+2-2+2-1+1-1+1-1+0+0-1+0-1+0-1+0-1+1-1+1+0+1+0+1+0+1+0+1  $\\ 
\hline 
$\mathbf{k}7_{16} $ &(-76, -41)& $-1+2-1+1-1+0+0+0+0+0+0+0+0+0+0+0-1+0-1+0-1+0-1+1-1+1-1+1+0+1+0+1+0+1+0+1  $\\ 
\hline 
$\mathbf{k}7_{17} $ &(44, 81)& $1+0+1+0+1+0+1+0+1-1+1-1+1-1+0-1+0-1+0-1+0+0+0+0+0+1-1+1-1+1-2+2-2+2-2+2-2+1  $\\ 
\hline 
$\mathbf{k}7_{18} $ &(-6, 3)& $1-1+1-1+0+0+0+1-1+1  $\\ 
\hline 
$\mathbf{k}7_{19} $ & (-7, 3)& $-1+1-1+1-1+2-1+2-1+1-1  $\\ 
\hline 
$\mathbf{k}7_{20} $ &(14, 28)& $1+0+1+0+1+0+0+0-1+0-1+1-1+1-1  $\\ 
\hline 
$\mathbf{k}7_{21} $ &(-49, -25)& $-1+1+0+0+0+0+0+0+0+0-1+0-1+0-1+0+0+0+1+0+1+0+1+0+1  $\\ 
\hline 
$\mathbf{k}7_{22} $ &(33, 63)& $1+0+1+0+1+0+0+0+0+0+0+0+0+0+0+0+0-0+0-1+0-1+0+0+0+0+0+0+0-1+1  $\\ 
\hline 
$\mathbf{k}7_{23} $ &(33, 64)& $1+0+1+0+1+0+0+0+0+0+0+0+0+0+0+0+0+0-1+0-1+0+0+0+0+0+0+0+0+0+1-1  $\\ 
\hline 
$\mathbf{k}7_{24} $ &(1, 12)& $1-1+1+0+0+0+0+0+1-1+1-1  $\\ 
\hline 
$\mathbf{k}7_{25} $ &(-3, 3)& $1-1+1-1+1-1+1  $\\ 
\hline 
$\mathbf{k}7_{26} $ &(-5, 6)& $-1+2-1+2-2+2-2+1-1+1-1+1  $\\ 
\hline 
$\mathbf{k}7_{27} $ &(23, 49)& $1+0+1+0+1-1+1-1+1-1+1-1+1-1+1-1+0+0-1+1-2+2-2+2-2+2-1  $\\ 
\hline 
$\mathbf{k}7_{28} $ &(-38, -22)& $-1+0-1+1-1+1-1+0+0+0+0+0+1+0+1+0+1  $\\ 
\hline 
$\mathbf{k}7_{29} $ & & Not computed  \\ 
\hline 
$\mathbf{k}7_{30} $ &(28, 59)& $1+0+1+0+1-1+1-1+1-1+1-1+1-1+1-1+1-1+1-1+0+0-1+1-2+2-2+2-2+2-2+1  $\\ 
\hline 
$\mathbf{k}7_{31} $ & & Not computed  \\ 
\hline 
$\mathbf{k}7_{32} $ &(-64, -34)& $-1+2-1+1-1+0+0+0+0+0+0+0+0-1+0-1+0-1+0-1+1-1+1+0+1+0+1+0+1+0+1  $\\ 
\hline 
$\mathbf{k}7_{33} $ &(-60, -32)& $1-1+0+0+0+0+0+0+0+0+0-1+0-1+0-1+0+0+0+0+0+0+1+0+1+0+1+0+1  $\\ 
\hline 
$\mathbf{k}7_{34} $ &(19, 38)& $1+0+1+0+1+0+0+0+0+0-1+1-2+1-1+0+0+0+1-1  $\\ 
\hline 
$\mathbf{k}7_{35} $ &(-62, -33)& $-1+2-1+1-1+0+0-1+1-1+1-1+1-2+1-1+0-1+0+0+0+0+0+1+0+1+0+1+0+1  $\\ 
\hline 
$\mathbf{k}7_{36} $ &(39, 73)& $1+0+1+0+1+0+1+0+0+0+0+0+0+0-1+1-2+1-2+1-1+0+0+0+0-1+2-2+2-1+0+0-1+2-1  $\\ 
\hline 
$\mathbf{k}7_{37} $ &(-53, -27)& $1-1+0+0+0+0+0-1+1-2+1-1+0+0+0+0+0+0+0+0+0+0 +1+0+1+0+1  $\\ 
\hline 
$\mathbf{k}7_{38} $ &(-80, -43)& $1-2+1-1+1+0+0+0+0+0+0+0+0+0+0+0+0+0-1+0-1+0-1+0-1+1-1+1-1+1+0+1+0+1+0+1+0+1  $\\ 
\hline 
$\mathbf{k}7_{39} $ &(8, 22)& $1+0+1+0+0+0-1+1-1+1-1+1-1+1-1  $\\ 
\hline 
$\mathbf{k}7_{40} $ &(-3, 7)& $-1+1-1+2-1+2-1+1-1+1-1  $\\ 
\hline 
$\mathbf{k}7_{41} $ &(-5, 6)& $1-1+1-1+1-2+2-2+3-2+2-1  $\\ 
\hline 
$\mathbf{k}7_{42} $ &(49, 90)& $1+0+1+0+1+0+1+0+0+0+0+0+0+0+0+0+0+0-1+0-1+0-1+0+0+0+0+0+0+0+0+0+1-1+1-1+0+0+0+0-1+1  $\\ 
\hline 
$\mathbf{k}7_{43} $ &(50, 91)& $1+0+1+0+1+0+1+0+0+0+0+0+0+0+0+0+0-0+0-1+0-1+0-1+0+0+0+0+0+0+0-1+1-1+1+0+0+0+0+0+1-1  $\\ 
\hline 
$\mathbf{k}7_{44} $ &(-4, 8)& $-1+2-3+3-2+3-2+2-1+1-1+1-1  $\\ 
\hline 
$\mathbf{k}7_{45} $ &(1, 11)& $1-1+2-2+3-3+3-3+2-2+1  $\\ 
\hline 
$\mathbf{k}7_{46} $ &(4, 15)& $1-1+2-2+3-3+3-3+3-2+1-1  $\\ 
\hline 
$\mathbf{k}7_{47} $ &(-39, -21)& $1-1+1-2+1-2+1-2+1-1+1-1+1+0+1+0+1+0+1  $\\ 
\hline 
$\mathbf{k}7_{48} $ &(24, 48)& $1+0+1+0+1+0+1+0+0+0+0-1+0-1+0-1+0-1+0+0+0+1+0+1-1  $\\ 
\hline 
$\mathbf{k}7_{49} $ &(-6, 4)& $-1+0+0+1+0+1+0+1-1+1-1  $\\ 
\hline 
$\mathbf{k}7_{50} $ &(22, 44)& $1+0+1+0+1+0+1-1+1-1+0-1+0+0-1+1-1+1-1+0+0-1+1  $\\ 
\hline 
$\mathbf{k}7_{51} $ &(-49, -25)& $-1+2-1+1-1+1-1+0+0-1+0-1+0-1+0+0+0+0+1+0+1+0 +1+0+1  $\\ 
\hline 
$\mathbf{k}7_{52} $ &(-38, -19)& $-1+1-1+2-2+2-2+1-2+1-2+1-1+1+0+1+0+1+0+1  $\\ 
\hline 
$\mathbf{k}7_{53} $ &(-30, -15)& $-1+0+0+0+0+0+0+0-1+1-1+1+0+1+0+1  $\\ 
\hline 
$\mathbf{k}7_{54} $ &(-2, 9)& $1-1+1-1+1-1+2-2+2-2+2-1  $\\ 
\hline 
$\mathbf{k}7_{55} $ &(-50, -28)& $1-1+1-1+0+0-1+0-1+0-1+0+0+0+0+0+1+0+1+0+1-0+1  $\\ 
\hline 
$\mathbf{k}7_{56} $ &(21, 40)& $1+0+1+0+1-1+1-1+1-1+1-1+1-1+1-2+1-1+1-1  $\\ 
\hline 
$\mathbf{k}7_{57} $ &(24, 44)& $1+0+1+0+1-1+1-1+1-1+1-1+1-1+1-1+0+0+0+0-1  $\\ 
\hline 
$\mathbf{k}7_{58} $ &(-61, -33)& $1-1+0-1+0+0+0+1+0+0+0+0-1+0-1+0-1+0+0+0+0+0+1+0+1+0+1+0+1  $\\ 
\hline 
$\mathbf{k}7_{59} $ &(32, 61)& $1+0+1+0+1+0+1+0+0+0+0+0-1+0-1+0-1+0+0+0+0+1-1+1-1+0+0-1+2-1  $\\ 
\hline 
$\mathbf{k}7_{60} $ &(-65, -35)& $-1+1+0+1-1+1-1+0+0+0-1+1-1+0-1+0+0-1+0+0+0+0+0+0+1+0+1+0+1+0+1  $\\ 
\hline 
$\mathbf{k}7_{61} $ &(29, 52) & $1+0+1+0+1+0+1+0+0+0+0-1+0-1+0-1+0+0+0+0+1-1+1-1 $\\
\hline 
$\mathbf{k}7_{62} $ &(37, 69)& $1+0+1+0+1+0+1+0+1-1+1-1+0-1+0-1+0-1+0+0+0+0+0+1-1+1-1+1-1+0+0-1+1  $\\ 
\hline 
$\mathbf{k}7_{63} $ &(-55, -28)& $-1+2-2+1+0+0+1-2+2-2+1-1+0-1+0+0-1+1-1+1-1+1+0+1+0+1+0+1  $\\ 
\hline 
$\mathbf{k}7_{64} $ &(38, 73)& $1+0+1+0+1+0+0+0+0+0+0+0+0+0+0+0-0+0+0+0-1+1-2+1-1+0+0+0+0+0+0+0+0-1+2-1  $\\ 
\hline 
$\mathbf{k}7_{65} $ &(-2, 5)& $1+0+0+0-1+1-1+1  $\\ 
\hline 
$\mathbf{k}7_{66} $ &(-60, -31)& $1-2+2-2+1-1+1+0+0+0+0-1+0+0-1+0+0+0-1+1-1+1-1+1+0+1+0+1+0+1  $\\ 
\hline 
$\mathbf{k}7_{67} $ &(-90, -50)& $1-1+0-1+1-1+1+0+0+0+1-1+1-1+1-1+0+0+0-1+0-1+0-1+0-1+0+0+0+0+0+0+1+0+1+0+1+0+1+0+1  $\\ 
\hline 
$\mathbf{k}7_{68} $ &(-28, -12)& $-1+1-1+1-1+0+0+0+0+0+0+0+0+0+1+0+1  $\\ 
\hline 
$\mathbf{k}7_{69} $ &(-75, -40)& $-1+2-1+0+0-1+2-2+2-1+1-1+1-1+0+0-1+0-1+0-1+0-1+1-1+1-1+1+0+1+0+1+0+1+0+1  $\\ 
\hline 
$\mathbf{k}7_{70} $ &(-1, 7)& $1-1+1-1+1+0+0+1-1  $\\ 
\hline 
$\mathbf{k}7_{71} $ &(-78, -42)& $1-2+1+0+0+0-1+1+0+0+1-1+0+0+0-1+0+0-1+0-1+1-1+0+0+0+0+0+0+0+1+0+1+0+1+0+1  $\\ 
\hline 
$\mathbf{k}7_{72} $ &(2, 13)& $1+0+0+0+0+0+0+0+1-1+1-1  $\\ 
\hline 
$\mathbf{k}7_{73} $ &(-88,  -49)& $1-1+0-1+1-1+1+0+0+0+1-1+0+0+0+0-1+1-1+0-1+0-1+1-2+1-1+1-1+1-1+1+0+1+0+1+0+1+0+1  $\\ 
\hline 
$\mathbf{k}7_{74} $ &(-1, 11)& $-1+3-3+3-3+3-2+2-1+1-1+1-1  $\\ 
\hline 
$\mathbf{k}7_{75} $ &(-33, -17)& $-1+1-1+1-1+0+0-1+0+0+0+0+1+0+1+0+1  $\\ 
\hline 
$\mathbf{k}7_{76} $ &(39, 75)& $1+0+1+0+1+0+1-1+1-1+1-1+1-1+1-1+1-2+1-1+0-1+1-1+0+0+0+0+0+0-1+1-1+2-2+2-1  $\\ 
\hline 
$\mathbf{k}7_{77} $ &(-10, 1)& $-1+1-1+1-1+2-1+2-2+1-1+1  $\\ 
\hline 
$\mathbf{k}7_{78} $ &(42, 80)& $1+0+1+0+1+0+1-1+1-1+1-1+1-1+1-1+1-1+0+0-1+0-1+1-2+2-2+2-2+2-2+1+0+0+1-2+2-2+1  $\\ 
\hline 
$\mathbf{k}7_{79} $ & & Not computed  \\ 
\hline 
$\mathbf{k}7_{80} $ &(3, 11)& $1+0+0+1-1+1-1+1-1  $\\ 
\hline 
$\mathbf{k}7_{81} $ &(10, 24)& $1+0+1+0+0+0+0+0+0+0-1+1-1+1-1  $\\ 
\hline 
$\mathbf{k}7_{82} $ &(-59, -30)& $1-2+1+0+0+0+0+0+0-1+0+0-1+1-1+0+0+0+0+0+0+0+0+ 0+0+1+0+1+0+1  $\\ 
\hline 
$\mathbf{k}7_{83} $ &(-23, -10)& $1-1+1-2+0-1+0+0+0+1+0+1+0+1  $\\ 
\hline 
$\mathbf{k}7_{84} $ &(-12, 0)& $1-1+1-1+1-2+2-3+2-2+3-1+1  $\\ 
\hline 
$\mathbf{k}7_{85} $ &(-2, 7)& $1-1+2-1+0+0-1+1-1+1  $\\ 
\hline 
$\mathbf{k}7_{86} $ &(1, 10)& $1-2+3-3+4-3+3-2+1-1  $\\ 
\hline 
$\mathbf{k}7_{87} $ &(-37, -19)& $1-1+0-1+0+0-1+1-1+0+0+0+0+0+1+0+1+0+1  $\\ 
\hline 
$\mathbf{k}7_{88} $ &(36, 69)& $1+0+1+0+1+0+0+0+0+0+0+0+0+0+0+0-0+0+0-1+0+0-1+1-1+0+0+0+0+0+0+1-2+1  $\\ 
\hline 
$\mathbf{k}7_{89} $ &(-4, 6)& $1-1+2-3+4-4+3-3+2-1+1  $\\ 
\hline 
$\mathbf{k}7_{90} $ &(-25, -10)& $-1+1-1+1-1+1-1+0+0+0+0+0+0+1+0+1  $\\ 
\hline 
$\mathbf{k}7_{91} $ &(-5, 5)& $1-1+2-3+3-4+4-3+3-2+1  $\\ 
\hline 
$\mathbf{k}7_{92} $ &(-1, 8)& $-1+2-1+1+0+0+1-1+1-1  $\\ 
\hline 
$\mathbf{k}7_{93} $ &(2, 13)& $1-1+2-3+4-4+4-3+3-2+1-1  $\\ 
\hline 
$\mathbf{k}7_{94} $ &(0, 9)& $1+0+0-1+1-1+2-1+1-1  $\\ 
\hline 
$\mathbf{k}7_{95} $ &(-10, -3)& $-1+1-2+2-1+2-1+1  $\\ 
\hline 
$\mathbf{k}7_{96} $ &(-6, -1)& $-1+1-1+2-1+1  $\\ 
\hline 
$\mathbf{k}7_{97} $ &(-7, 0)& $-1+2-2+2-2+2-1+1  $\\ 
\hline 
$\mathbf{k}7_{98} $ &(-12, -4)& $-1+0+0+0+1-1+2-1+1  $\\ 
\hline 
$\mathbf{k}7_{99} $ &(-14, -5)& $1-2+1-2+1+0+0+1+0+1  $\\ 
\hline 
$\mathbf{k}7_{100} $ &(-12, -1)& $-1+1-1+1+0+1-1+1-2+2-2+2  $\\ 
\hline 
$\mathbf{k}7_{101} $ &(27, 52)& $1+0+1+0+1+0+1+0+0+0+0+0-1+0-1+0-1+0-1+0+0+0+1+0+1-1  $\\ 
\hline 
$\mathbf{k}7_{102} $ &(8, 20)& $1+0+1+0+1-1+1-2+2-3+2-2+1  $\\ 
\hline 
$\mathbf{k}7_{103} $ &(-8, 2)& $1-2+3-4+4-4+4-3+2-1+1  $\\ 
\hline 
$\mathbf{k}7_{104} $ &(0, 10)& $1+0+0+0+1-1+1-2+1-1+1  $\\ 
\hline 
$\mathbf{k}7_{105} $ &(-3, 8)& $-1+2-1+2-1+1-1+0-1+1-1+1  $\\ 
\hline 
$\mathbf{k}7_{106} $ &(8, 21)& $1+0+1+0+1-1+1-1+0-1+0+0-1+1  $\\ 
\hline 
$\mathbf{k}7_{107} $ &(-10, 2)& $-1+1-2+1+0+1+0+1+0+1-1+1-1  $\\ 
\hline 
$\mathbf{k}7_{108} $ &(-14, -3)& $-1+1-2+3-4+5-4+4-3+2-1+1  $\\ 
\hline 
$\mathbf{k}7_{109} $ &(-47, -23)& $-1+1+1-1+1-1+1-2+1-1+0-1+1-2+1-1+1-1+1+0+1+0+1+0+1  $\\ 
\hline 
$\mathbf{k}7_{110} $ &(12, 28)& $1+0+1+0+1-1+1-1+0+0-1+1-2+2-1+1-1  $\\ 
\hline 
$\mathbf{k}7_{111} $ &(-57, -30)& $1-1+0-1+0+0+0+1+0+0+0+0-1+0-1+0-1+0+0+0+0+1+0+1+0+1+0+1  $\\ 
\hline 
$\mathbf{k}7_{112} $ &(26, 51)& $1+0+1+0+1+0+1-1+1-1+1-2+1-1+0-1+1-1+0+1-2+2-2+3-3+1  $\\ 
\hline 
$\mathbf{k}7_{113} $ &(7, 18)& $1+0+1+0+1-2+2-3+3-3+2-1  $\\ 
\hline 
$\mathbf{k}7_{114} $ &(0, 9)& $1+0+1+0-1+0-1+1-1+1  $\\ 
\hline 
$\mathbf{k}7_{115} $ &(14, 31)& $1+0+1+0+1-1+1-1+0+0+0+0+0-1+1-2+2-1  $\\ 
\hline 
$\mathbf{k}7_{116} $ &(-10, 1)& $-1+1-1+1+0+1-1+1-1+1-1+1  $\\ 
\hline 
$\mathbf{k}7_{117} $ &(-3, 6)& $-1+1-1+2+0+0+1-1+1-1  $\\ 
\hline 
$\mathbf{k}7_{118} $ &(7, 18)& $1+0+1+0+1-1+1-3+2-2+2-1  $\\ 
\hline 
$\mathbf{k}7_{119} $ &(15, 34)& $1+0+1+0+1-1+1-1+1-1+1-1+0+0-1+1-2+1-1+1  $\\ 
\hline 
$\mathbf{k}7_{120} $ &(2, 9)& $1-1+3-3+3-3+2-1  $\\ 
\hline 
$\mathbf{k}7_{121} $ &(-6, -1)& $-1+1-1+2-1+1  $\\ 
\hline 
$\mathbf{k}7_{122} $ &(-32, -15)& $1-2+2-2+1-1+0+0-1+0+0+0+0+1+0+1+0   $\\ 
\hline 
$\mathbf{k}7_{123} $ &(44, 84)& $1+0+1+0+1+0+1-1+1-1+1-1+1-1+1-1+1-1+0+0+0-1+0+0-1+0+0+0+0+0+0-1+1+0+0+0-1+2-2+2-1  $\\ 
\hline 
$\mathbf{k}7_{124} $ &  (8, 19)& $1+0+1+0+1-2+2-3+3-3+3-2 $\\ 
\hline 
$\mathbf{k}7_{125} $ &(0, 9)& $2-1+1-1+0+0-1+1-1+1  $\\ 
\hline 
$\mathbf{k}7_{126} $ &(-41, -19)& $1-2+2-2+1-1+0+0+0+0+0-1+1-1+1-1+1-1+1+0+1+0+1  $\\ 
\hline 
$\mathbf{k}7_{127} $ &(-7, -1)& $1-2+2-3+3-2+2  $\\ 
\hline 
$\mathbf{k}7_{128} $ &(-5, 4)& $1-1+1-1+0+1-1+1-1+1  $\\ 
\hline 
$\mathbf{k}7_{129} $ &(-4, 8)& $-1+2-2+2-1+1-1+1+0+1-1+1-1  $\\ 
\hline

 \end{longtable}

\newpage

 \begin{longtable}{|l||l|l|p{3.5in}|}

 \caption{Dowker-Thistlethwaite codes for 57 census knots}\\
  \hline
   & OC & C & DT code \\
  \hline 
  \endfirsthead

 \caption[]{\emph{continued}}\\
 \hline
  & OC & C & DT code \\
   \hline 
 \endhead

 \hline
 \endlastfoot

 $k7_{1}$    &16 &     12 &     -4 -14 -24 -22 -20 -18 -16 -2 -12 -10 -8 -6  \\ \hline
$k7_{2}$    &25 &     13 &     -4 -16 -26 -24 -22 -20 -18 -2 -14 -12 -10 -8 -6  \\ \hline
$k7_{10}$    &321&     22 &     4 18 26 -24 40 38 34 30 2 -28 -32 -36 42 -44 16 -20 14 -22 12 10 8 -6 \\ \hline
$k7_{12}$    &329 &      23 &     4 18 26 -42 24 -40 -38 -34 -30 2 32 36 8 -44 -46 20 -16 22 -14 -12 -10 -6 -28 \\ \hline
$k7_{18}$    &521 &    15 &    -4 -12 -16 24 -14 22 -2 -8 26 28 30 10 6 18 20  \\ \hline
$k7_{24}$    &595 &    18 &    4 16 22 -20 30 -18 28 2 -26 -10 32 -34 -36 14 12 8 -6 -24 \\ \hline
$k7_{25}$    &600 &    16 &    4 10 12 18 -16 2 -28 -32 -24 26 30 8 -6 20 -14 22  \\ \hline
$k7_{26}$    &656 &    22 &    4 14 20 28 -26 36 -22 2 -34 24 -44 16 -10 -38 40 42 12 -18 8 -6 30 32 \\ \hline
$k7_{40}$    &960 &    12 &    -4 -8 -16 -2 18 20 22 -6 24 10 12 14  \\ \hline
$k7_{41}$    &1063 &    18 &    4 10 -14 30 2 -22 -8 28 -24 32 -12 -18 -34 -36 6 20 16 -26 \\ \hline
$k7_{44}$    &1126 &    21 &    4 14 20 -18 -26 34 2 -24 30 38 -40 -42 -32 36 -8 16 12 10 28 -6 -22 \\ \hline
$k7_{45}$    &1217 &    10 &    4 12 14 16 18 20 2 6 8 10  \\ \hline
$k7_{46}$    &1243 &    11 &    10 14 16 18 20 22 2 4 6 8 12  \\ \hline
$k7_{51}$    &1392 &    62 &    -14 74 110 34 -80 -116 -40 -28 56 94 -8 -62 -100 2 72 108 -20 -78 -114 26 102 64 120 -84 -86 122 -106 -30 -4 112 36 10 -118 48 -52 -90 -16  58 96 22 -44 -46 -66 -50 124 54 -32 -6 -60 38 12 88 70 -92 -18 -76 98 24 -42 82 -68 -104  \\ \hline
$k7_{54}$    &1434 &18& 4 10 -16 34 -18 2 20 -22 -36 28 -30 32 6 8 -12 14 24 26 \\ \hline
$k7_{55}$    &1547 &    178 &    -4 42 -146 38 188 268 58 224 304 94 234 -68 -278 -206 -138 -348 182 110 320 -252 2 148 -150 152 -154 156 -158 86 296 118 328 200 130 340 -174 -102 -312 -244 -32 288 216 -8 -160 -192 194 122 332 -64 -274 238 26 176 -210 -142 -352 -80 -290 324 -164 -14 -196 -228 230 -232 -20 -170 -98 -308 346 74 284 212 356 40 254 -256 258 -260 262 -264 190 -12 222 306 96 236 24 -280 -208 -140 -350 -78 112 322 -266 -56 -298 300 226 -18 -168 134 344 72 282 -316 -248 -36 -186 -114 -270 -60 -302 -92 -334 336 -338 -66 -276 -204 -136 30 180 108 318 6 -44 46 -48 50 -52 54 84 294 116 326 202 132 342 70 -104 -314 -246 -34 -184 218 -10 -162 88 -90 120 330 -62 -272 240 28 178 106 -144 -354 -82 -292 -220 -166 -16 -198 124 -126 128 -22 -172 -100 -310 -242 76 286 214 -250  \\ \hline
$k7_{59}$    &1690 &    114 &    -26 -196 -100 -32 202 106 38 -208 -112 -44 214 118 50 60 -164 -98 66 134 172 -72 -140 -178 78 146 184 -84 -152 -190 -96 -162 2 -130 -168 -8 136 174 14 -142 -180 -20 148 186 -220 -222 224 -226 228 192 126 166 64 132 -36 -70 -138 42 76 144 -48 -82 -150 54 -56 -30 4 198 -170 -10 -204 176 16 210 -182 -22 -216 188 86 88 -90 92 -94 194 128 6 200 104 -12 -206 -110 18 212 116 -24 -218 -122 124 -58 -28 62 -102 -34 -68 108 40 74 -114 -46 -80 120 52 154 -156 158 -160  \\ \hline
$k7_{65}$    &1811 &    17 &    4 12 16 -26 14 -24 2 8 -28 -30 -32 -34 -10 -6 -18 -20 -22 \\ \hline
$k7_{67}$    &1839 &    204 &    -34 -330 -234 -180 -86 338 242 188 94 -346 -250 -196 -102 354 258 204 110 126 312 -222 132 318 -228 2 -142 -274 -370 -10 150 282 378 18 -158 -290 -386 -26 166 298 394 124 -364 -176 -270 54 148 280 376 -62 -156 -288 -384 70 164 296 392 -78 -172 -304 -400 -268 -174 -362 408 220 38 -316 226 44 -322 -326 4 332 -276 -372 -12 -340 284 380 20 348 -292 -388 -28 -356 300 396 404 216 366 -272 52 146 278 -92 -60 -154 -286 100 68 162 294 -108 -76 -170 -302 116 -266 -120 -82 -36 314 130 -42 320 136 324 328 50 144 -184 -90 -58 -152 192 98 66 160 -200 -106 -74 -168 208 114 402 308 368 6 334 238 -374 -14 -342 -246 382 22 350 254 -390 -30 -358 -262 398 -212 -118 -84 -178 -128 40 224 -134 46 230 -138 232 -140 48 -236 -182 -88 -56 244 190 96 64 -252 -198 -104 -72 260 206 112 406 218 310 8 336 240 186 -16 -344 -248 -194 24 352 256 202 -32 -360 -264 -210 306 214 -122 -80  \\ \hline
$k7_{70}$    &1935 &    13 &    -4 -8 -10 18 -2 20 22 24 26 6 12 14 16  \\ \hline
$k7_{72}$    &1964 &    16 &    4 16 22 -20 28 -18 26 2 -24 -10 30 -32 14 12 8 -6  \\ \hline
$k7_{74}$    &1971 &    21 &    4 16 20 32 -18 28 26 2 -24 30 -34 36 14 12 10 -38 40 -42 8 -6 22 \\ \hline
$k7_{77}$    &2001 &    17 &    4 10 -14 20 2 -28 -8 22 -24 26 6 -34 32 -30 -12 -18 16 \\ \hline
$k7_{80}$    &2166 &    10 &    4 12 -16 14 -18 2 8 -20 -10 -6  \\ \hline
$k7_{84}$    &2257 &    22 &    4 10 -18 28 2 -16 24 -36 -8 30 -32 34 -38 12 6 -44 42 -40 -26 -14 -22 20  \\ \hline
$k7_{85}$    &2272 &    13 &    -4 -8 12 -2 20 6 26 24 22 10 18 16 14  \\ \hline
$k7_{86}$    &2284 &    9 &    6 12 14 18 16 4 2 10 8  \\ \hline
$k7_{89}$    &2362 &    10 &    6 14 12 20 18 16 4 2 10 8  \\ \hline
$k7_{91}$    &2488 &    10 &    -6 -12 -14 -20 -18 -16 -4 -2 -10 -8  \\ \hline
$k7_{92}$    &2508 &    14 &    4 12 16 -14 22 2 -20 24 -26 -28 10 8 -6 -18  \\ \hline
$k7_{93}$    &2520 &    11 &    6 14 16 22 20 18 2 4 12 10 8  \\ \hline
$k7_{94}$    &2543 &    12 &    4 12 -16 -20 -18 2 22 24 -10 -8 -6 14  \\ \hline
$k7_{95}$    &2553 &    10 &    -4 -8 14 -2 16 18 20 6 12 10  \\ \hline
$k7_{96}$    &2552 &    11 &    -4 -8 14 -2 16 18 20 6 10 22 12  \\ \hline
$k7_{97}$    &2623 &    9 &    -4 -8 -10 -14 -2 16 -6 18 12  \\ \hline
$k7_{98}$    &2624 &    12 &    -4 -8 16 -2 18 20 22 24 6 10 14 12  \\ \hline
$k7_{99}$    &2642 &    12 &    -6 -10 -12 18 -2 -4 20 22 24 8 14 16  \\ \hline
$k7_{100}$    &2743 &    19 &    4 10 -20 28 2 -18 26 22 -36 -8 30 -34 14 12 6 -38 -16 -24 -32 \\ \hline
$k7_{103}$    &2858 &    10 &    -6 -14 -12 -16 -18 -20 -4 -2 -8 -10  \\ \hline
$k7_{104}$    &2861 &    13 &    4 14 18 16 -24 -22 -20 2 6 -26 -12 -10 -8  \\ \hline
$k7_{105}$    &2869 &    18 &    4 14 -20 18 24 -32 22 2 30 -28 -16 10 34 36 6 -12 -8 26 \\ \hline
$k7_{107}$    &2888 &    21 &    4 10 -18 20 2 22 -24 38 42 30 -32 34 -36 -40 6 -8 -12 14 -26 16 -28 \\ \hline
$k7_{108}$    &2894 &    11 &    -8 -14 -16 -18 -22 -20 -2 -4 -6 -12 -10  \\ \hline
$k7_{109}$    &2900 &    53 &    -12 90 50 -96 -76 44 68 36 -6 -94 -54 10 -78 60 -102 -84 66 -48 -4 -92 20 -22 106 88 -16 -72 -40 42 8 -80 -100 30 -86 -46 -2 38 18 -56 -24 98 -28 -62 104 -34 -14 -70 52 74 -58 -26 82 -64 -32  \\ \hline
$k7_{110}$    &2925 &    27 &    6 10 -38 -26 2 -28 30 -32 34 -36 -50 -52 -54 -8 40 -42 44 -46 48 -4 -12 14 -16 18 -20 -22 -24 \\ \hline
$k7_{114}$    &3093 &    16 &    -4 -10 14 -20 -2 26 8 32 28 -24 -6 -30 12 18 16 -22  \\ \hline
$k7_{115}$    &3105 &    31 &    4 10 -26 44 2 -28 -30 32 -34 36 -54 -58 -62 -8 -12 -46 48 -50 52 56 60 6 14 -16 18 -20 38 -22 40 -24 42 \\ \hline
$k7_{116}$    &3169 &    17 &    4 10 -18 24 2 -16 20 -32 -8 26 -30 12 6 -34 -14 -22 -28 \\ \hline
$k7_{117}$    &3195 &    12 &    4 10 12 -18 2 8 -20 -22 -24 -6 -14 -16  \\ \hline
$k7_{118}$    &3199 &    17 &    4 10 -16 18 2 -28 -30 -34 -24 26 32 6 -8 20 -12 -14 22 \\ \hline
$k7_{119}$    &3234 &    33 &    10 14 -18 20 56 -32 2 -34 36 50 -38 40 -42 44 -64 -66 -12 -48 52 -54 58 -60 62 16 4 -6 22 8 -24 26 -28 46 -30 \\ \hline
$k7_{120}$    &3310 &    7 &    4 10 12 14 2 8 6  \\ \hline
$k7_{121}$    &3320 &    12 &    -4 -10 16 18 -2 -20 -22 -24 8 6 -12 -14  \\ \hline
$k7_{125}$    &3423 &    15 &    4 14 20 -26 18 16 -24 2 10 8 -28 -30 -12 -6 -22  \\ \hline
$k7_{127}$    &3505 &    8 &    -4 -8 12 -2 -14 6 -16 -10  \\ \hline
$k7_{128}$    &3536 &    11 &    4 8 -12 -16 2 -18 -22 -20 -10 -6 -14  \\ \hline
$k7_{129}$    &3541 &    21 &    4 14 -18 22 -20 30 2 38 42 -12 32 -34 -36 -40 6 10 -8 -24 -26 16 -28 \\ \hline

 \end{longtable}




















\nonumsection{Acknowledgements}
This project could only be accomplished using computer tools developed over many years by Nathan Dunfield, Oliver Goodman, Morwen Thistlethwaite,
and Jeff Weeks.  We thank them not only for their programs, but for
readily and frequently assisting us with various aspects of this project.
We thank Walter Neumann for his encouragement and advice.  We thank M.
Ochiai for computing Jones polynomials for $\mathbf{k}7_{61}$ and $\mathbf{k}7_{106}$. We thank Morwen Thistlethwaite for providing link pictures. We thank 
Rob Scharein for making KnotPlot pictures of many of these new knots 
\cite{scharein}.  This
project began as an NSF VIGRE summer project at Columbia University, and
we appreciate NSF support, and thank the participants: Gerald Brant,
Jonathan Levine, Rustam Salari.

\nonumsection{References}
\bibliography{main.bib}

\begin{thebibliography}{10}

\bibitem{cdw99}
Patrick~J. Callahan, John~C. Dean, and Jeffrey~R. Weeks.
\newblock The simplest hyperbolic knots.
\newblock {\em J. Knot Theory Ramifications}, 8(3):279--297, 1999.

\bibitem{dunfield}
Nathan Dunfield.
\newblock Additions to snappea.\\
\newblock {\em http://abel.math.harvard.edu/$\sim$nathand/snappea/index.html}.

\bibitem{gs99}
Robert~E. Gompf and Andr{\'a}s~I. Stipsicz.
\newblock {\em {$4$}-manifolds and {K}irby calculus}, volume~20 of {\em
  Graduate Studies in Mathematics}.
\newblock American Mathematical Society, Providence, RI, 1999.

\bibitem{snap}
Oliver Goodman.
\newblock Snap.
\newblock {\em http://www.ms.unimelb.edu.au/$\sim$snap/}.

\bibitem{hw89}
Martin Hildebrand and Jeffrey Weeks.
\newblock A computer generated census of cusped hyperbolic {$3$}-manifolds.
\newblock In {\em Computers and mathematics (Cambridge, MA, 1989)}, pages
  53--59. Springer, New York, 1989.

\bibitem{testisom}
Derek Holt and Sarah Rees.
\newblock isom\_quotpic.\\
\newblock {\em ftp://ftp.maths.warwick.ac.uk/people/dfh/isom\_quotpic/}.

\bibitem{Knotscape}
Jim Hoste and Morwen Thistlethwaite.
\newblock Knotscape 1.01.\\
\newblock {\em http://dowker.math.utk.edu/knotscape.html}.

\bibitem{lackenby}
Marc Lackenby.
\newblock The volume of hyperbolic alternating link complements.
\newblock {\em arXiv:math.GT/0012185}.

\bibitem{lackenby2000}
Marc Lackenby.
\newblock Word hyperbolic {D}ehn surgery.
\newblock {\em Invent. Math.}, 140(2):243--282, 2000.

\bibitem{meyerhoff92}
Robert Meyerhoff.
\newblock Geometric invariants for {$3$}-manifolds.
\newblock {\em Math. Intelligencer}, 14(1):37--53, 1992.

\bibitem{mm2001}
Hitoshi Murakami and Jun Murakami.
\newblock The colored {J}ones polynomials and the simplicial volume of a knot.
\newblock {\em Acta Math.}, 186(1):85--104, 2001.

\bibitem{ochiai}
M~Ochiai.
\newblock K2k version 1\_3\_1.\\
\newblock {\em http://amadeus.ics.nara-wu.ac.jp/$\sim$ochiai/download.html}.

\bibitem{rolfsen}
Dale Rolfsen.
\newblock {\em Knots and links}, volume~7 of {\em Mathematics Lecture Series}.
\newblock Publish or Perish Inc., Houston, TX, 1990.
\newblock Corrected reprint of the 1976 original.

\bibitem{saveliev99}
Nikolai Saveliev.
\newblock {\em Lectures on the topology of {$3$}-manifolds}.
\newblock de Gruyter Textbook. Walter de Gruyter \& Co., Berlin, 1999.
\newblock An introduction to the Casson invariant.

\bibitem{scharein}
Rob Scharein.
\newblock {\em http://www.cecm.sfu.ca/$\sim$scharein/nhypcen/}.

\bibitem{snappea}
Jeff Weeks.
\newblock Snappea.
\newblock {\em http://www.northnet.org/weeks}.

\end{thebibliography}
\bibliographystyle{plain}

\end{document}